\documentclass[12pt,a4paper]{article}
\usepackage{amsmath,amsfonts,amssymb,amscd}

\def\gr{\operatorname{gr}}
\def\Ker{\operatorname{Ker}}

\def\id{\operatorname{id}}

\def\red{\operatorname{red}}

\def\GL{\operatorname{GL}}

\def\pr{\operatorname{pr}}
\def\Im{\operatorname{Im}}

\def\Int{\operatorname{Int}}

\def\mod{\operatorname{mod}}
\def\rank{\operatorname{rank}}

\newcounter{th}
\def\t{\refstepcounter{th}{\bf \noindent{Theorem} \arabic{th}. }}

\newcounter{le}
\def\l{\refstepcounter{le}{\bf \noindent{Proposition} \arabic{le}. }}

\newcounter{lem}
\def\lem{\refstepcounter{lem}{\bf \noindent{Lemma} \arabic{lem}. }}

\newcounter{de}

\newcounter{ex}

\begin{document}

\begin{center}
    \Large{\bf  Locally free sheaves on complex supermanifolds}\footnote{
    This work was partially supported by MPI Bonn, SFB TR $|$$12$, DFG 1388 and by the Russian Foundation
for Basic Research (grant no. 11-01-00465a). }
\end{center}

\begin{center}
   A.L.~Onishchik,  E.G.~Vishnyakova
\end{center}

\begin{center}
{\bf 1. Introduction}
\end{center}
\medskip

 An important part of the classical theory of real or complex
manifolds is the theory of (smooth, real analytic or complex
analytic) vector bundles. With any vector bundle over a manifold
$(M,\mathcal F)$ the sheaf of its (smooth, real analytic or complex
analytic) sections is associated which is a locally free sheaf of
$\mathcal F$-modules, and in this way all the locally free sheaves
of $\mathcal F$-modules over $(M,\mathcal F)$ can be obtained. In
the present paper, locally free sheaves of $\mathcal O$-modules on a
complex analytic supermanifold $(M,\mathcal O)$ (or equivalently
sheaves of sections of vector bundles over $(M,\mathcal O)$) are
studied.

It is well-known that any smooth supermanifold $(M,\mathcal{O})$ is
split, i.e. $\mathcal{O} \simeq \bigwedge_{\mathcal{F}}
\mathcal{G}$, where $\mathcal{G}$ is the sheaf of sections of a
certain vector bundle over $M$. In the complex case this statement
is false, see \cite{Green}. However, we can assign the split
supermanifold $(M,\gr\mathcal{O})$ to any complex analytic
supermanifold $(M,\mathcal{O})$, which is called the {\it retract}
of $(M,\mathcal{O})$. Given a locally free sheaf $\mathcal E$ of
$\mathcal O$-modules on a complex analytic supermanifold
$(M,\mathcal O)$, we construct a locally free sheaf $\gr
\mathcal{E}$ on the retract $(M,\gr\mathcal O)$, which is called the
{\it retract} of $\mathcal E$. It can be easily shown that $\gr
\mathcal{E} \simeq \gr\mathcal O \otimes \mathcal{E}_{\red}$, where
$\mathcal{E}_{\red}$ is the pullback of $\mathcal{E}$ with respect
to the natural embedding of the manifold $(M, \mathcal{F})$ into
$(M,\mathcal{O})$. In Section $2$ we obtained a classification of
locally free sheaves $\mathcal{E}$ of $\mathcal{O}$-modules which
have a given retract $\gr \mathcal{E}$ in terms of non-abelian
1-cohomology, Theorem \ref{Aut_E}. In the special case $\mathcal O
\simeq \gr\mathcal O$ our classification result can be simplified,
Theorem \ref{Aut_E_(gr)}.

In Section $3$ we study locally free sheaves of modules over
projective superspaces. In the case of complex projective spaces,
the problem of the
 (indecomposable) bundle classification is far from being solved, see
 \cite{Oko}.
 There are two cases, however, in which all bundles are known to be
 direct sums of line bundles~--- over $\mathbb{CP}^1$ by the classical
 Birkhoff -- Grothendieck Theorem and over $\mathbb{CP}^{\infty}$
 by the Barth -- Van de Ven -- Tyurin theorem.
 We study similar
question in the super context. In the case of $\mathbb{CP}^{1|m}$,
$m>0$, we showed that the Birkhoff -- Grothendieck Theorem does not
hold true. (The fact that this theorem is false for some
$\mathbb{CP}^{1|m}$ was noticed in \cite{Man}.) Furthermore, we
achieved the result similar to the Barth -- Van de Ven -- Tyurin
Theorem for projective superspaces.

Section $4$ is devoted to the study of the tangent sheaf
$\mathcal{T}$ of a split supermanifold $(M,\bigwedge \mathcal{G})$
in more details. The main result is here the equivalence of the
triviality of the $1$-cohomology class corresponding to
$\mathcal{T}$ and the existence of a holomorphic connection of the
bundle corresponding to the locally free sheaf of
$\mathcal{F}$-modules $\mathcal{G}$.

In Subsection $5$ a spectral sequence which connects the cohomology
with values in a locally free sheaf of $\mathcal{O}$-modules
$\mathcal{E}$ with the cohomology with values in its retract
$\gr\mathcal{E}$ is constructed. This spectral sequence permits to
compute the cohomology group $H^*(M,\mathcal{E})$ using the
cohomology class corresponding to $\mathcal{E}$ by Theorem
\ref{Aut_E_(gr)} and the cohomology group $H^*(M,\gr\mathcal{E})$.
Note that $\gr\mathcal{E}$ is a sheaf of sections of a certain
vector bundle over $M$. Hence to compute $H^*(M,\gr\mathcal{E})$ we
may use the well elaborated tools of complex analytic geometry. We
described the first two terms of the spectral sequence and the first
non zero differential.

A classification of locally free sheaves of $\mathcal{O}$-modules
over a smooth supermanifold $(M,\mathcal{O})$ was obtained in
\cite{Rempel}, Section $4.3$. It was shown that any locally free
sheaf of $\mathcal{O}$-modules $\mathcal{E}$ is isomorphic to $\gr
\mathcal{E}$. The similar result for fibre superbundles was proved
in \cite{Schmitt}. In \cite{Czyz} the split holomorphic case was
studied. In particular it was shown there that there exists a
holomorphic locally free sheaf of $\mathcal{O}$-modules over a
holomorphic supermanifold $(M,\mathcal{O})$, which is not isomorphic
to  its retract $\gr \mathcal{E}$. There a classification up to
isomorphism of locally free sheaves of $\mathcal{O}$-modules over a
(holomorphic) split supermanifold $(M,\mathcal{O})$,
$\mathcal{O}\simeq \bigwedge \mathcal{G}$, is obtained in terms of
cohomology set $H^1(M,\GL(n,\bigwedge \mathcal{G}))$. In the present
paper we suggest the different approach to the classification of
locally free sheaves of $\mathcal{O}$-modules over a split
supermanifold, Theorem \ref{Aut_E_(gr)}, and more generally over a
non-split supermanifold, Theorem \ref{Aut_E}. Let us explain the
difference in more details. Clearly one has a split homomorphism
$T:\GL(n,\bigwedge \mathcal{G})\to \GL(n,\mathbb{C})$ by taking the
degree zero part of $\GL(n,\bigwedge \mathcal{G})$. It induces the
map $H^1(T):H^1(M,\GL(n,\bigwedge \mathcal{G}))\to
H^1(M,\GL(n,\mathbb{C}))$. Denote by $a_{\mathcal{E}}$ the element
of $H^1(M,\GL(n,\bigwedge \mathcal{G}))$, which corresponds to a
locally free sheaf of $\mathcal{O}$-modules $\mathcal{E}$. Then, in
our notations, $\mathcal{E}_{\red}$ corresponds to
$H^1(T)(a_{\mathcal{E}})$. In our paper we classify all locally free
sheaves $\mathcal{E}$ such that $\mathcal{E}_{\red}$ is fixed.
Therefore, instead of computing $H^1(M,\GL(n,\bigwedge
\mathcal{G}))$, we suggest to use results concerning classification
of holomorphic bundles over a manifold, obtained in classical
geometry, and consider locally free sheaves with given retract on a
split supermanifold. The idea to classify non-split object, more
precisely, supermanifolds, using retrcts appeared firstly in
\cite{Green}.

We would like also to mention that, as in the classical case, the
line superbundles can be described using the $\exp$-map, see e.g.
\cite{BB}, Chapter VI, Section $2$. The Picard groups of generic
super-grassmannians were computed in \cite{Penkov_Skornyakov}.

\medskip

\noindent  {\bf Notations.}
\medskip

\begin{tabular}{cl}
  $(M,\mathcal{O})$ & supermanifold \\
   $(M,\gr\mathcal{O})$ & the retract of $(M,\mathcal{O})$ \\
   $\mathcal{T}=\mathcal{D}er \mathcal{O}$ & the tangent sheaf of $(M,\mathcal{O})$\\
    $\mathcal{A}ut \mathcal{O}$  & the sheaf of automorphisms of the structure sheaf $\mathcal{O}$\\
$\mathcal{A}ut_0\gr\mathcal{O}$  & the sheaf of automorphisms of $\gr\mathcal{O}$ preserving\\
&  the $\mathbb{Z}$-grading of $\gr \mathcal{O}$\\
$\gr \mathcal{E}$ & the retract of a locally free sheaf of $\mathcal{O}$-modules $\mathcal{E}$\\
$\mathcal{A}ut^{\mathcal{R}}\mathcal{E}$ &  the sheaf of automorphisms of a sheaf of $\mathcal{R}$-modules $\mathcal{E}$\\
$\mathcal{A}ut^{\mathcal{R}}_0\gr\mathcal{E}$ & the sheaf of automorphisms of a $\mathbb{Z}$-graded sheaf of   \\
&$\mathcal{R}$-modules $\gr\mathcal{E}$ preserving the $\mathbb{Z}$-grading of $\gr \mathcal{E}$\\
 $\mathcal{QA}ut \mathcal{E}$  & the sheaf of quasi-automorphisms of a locally free sheaf \\
 &of $\mathcal{O}$-modules $\mathcal{E}$\\
$\mathcal{QA}ut_0\gr\mathcal{E}$ & the sheaf of quasi-automorphisms of a $\mathbb{Z}$-graded locally  \\
&free sheaf $\gr\mathcal{E}$ preserving the $\mathbb{Z}$-grading of $\gr \mathcal{E}$\\
$\mathcal{A}ut^{\mathcal{F}}_{\bar 0}\mathcal{S}$ & a subsheaf of
$\mathcal{A}ut^{\mathcal{F}}\mathcal{S}$ consisting of even automorphisms \\
& of a $\mathbb{Z}_2$-graded sheaf $S$\\
$\mathcal{E}nd^{\mathcal{O}}\mathcal{E}$ & the sheaf of endomorphisms of a sheaf of $\mathcal{O}$-modules $\mathcal{E}$\\
\end{tabular}

\bigskip

\noindent  {\bf Acknowledgment.} The idea to study locally free sheaves of modules over complex supermanifolds was
inspired by communications of the second author with I.B.~Penkov during the Summer School
"Structures in Lie Representation Theory", Bremen 2009, Germany.
The authors are also grateful to V.~Serganova for useful discussions.

\bigskip

\begin{center}
{\bf 2. Main definitions and classification theorems}
\end{center}

\medskip

\noindent {\it 2.1 Main definitions and classification of complex
supermanifolds with a given retract}

We consider here complex analytic supermanifolds in the sense of
Berezin and Leites (see \cite{BL, ley}). Thus, a {\it supermanifold}
$(M,\mathcal{O})$ of dimension $n|m$ is a $\mathbb{Z}_2$-graded
ringed space which is locally isomorphic to a superdomain in
$\mathbb{C}^{n|m}$. The underlying complex manifold
$(M,\mathcal{F})$ is called the {\it reduction} of
$(M,\mathcal{O})$. Sometime we will denote it by $M$.  A {\it
morphism} $(M,\mathcal{O}_M) \to (N,\mathcal{O}_N)$ between two
supermanifolds with reductions $(M,\mathcal{F}_M)$ and
$(N,\mathcal{F}_N)$ is a morphism between $\mathbb{Z}_2$-graded
ringed spaces, i.e., a pair $F = (F_{red},F^*)$, where $F_{red}:
M\to N$ is a continuous mapping and $F^*: \mathcal{O}_N\to
(F_{red})_*\mathcal{O}_M$ is a homomorphism of sheaves of
$\mathbb{Z}_2$-graded ringed spaces. A morphism $F$ is called an
{\it isomorphism} if $F$ is invertible.

We consider $\mathbb{Z}_2$-graded sheaves of $\mathcal{O}$-modules
$\mathcal{S} = \mathcal{S}_{\bar 0}+\mathcal{S}_{\bar 1}$ on
$(M,\mathcal{O})$. Denote by $\Pi(\mathcal{S})$ the same sheaf of
$\mathcal{O}$-modules $\mathcal{S}$ supplied with the following
$\mathbb{Z}_2$-grading:
$$
\Pi(\mathcal{S})_{\bar 0} = \mathcal{S}_{\bar 1},\;
\Pi(\mathcal{S})_{\bar 1} = \mathcal{S}_{\bar 0}.
$$

A $\mathbb{Z}_2$-graded sheaf of $\mathcal{O}$-modules on
$(M,\mathcal{O})$ is called {\it free} ({\it locally free}) {\it of
rank} $p|q,\, p,q\ge 0$, if it is isomorphic (respectively, locally
isomorphic) to the $\mathbb{Z}_2$-graded sheaf of
$\mathcal{O}$-modules $\mathcal{O}^p\oplus\Pi(\mathcal{O})^q$. For
example, the tangent sheaf $\mathcal{T}$ of a supermanifold
$(M,\mathcal{O})$ of dimension $n|m$ is a locally free sheaf of
$\mathcal{O}$-modules of rank $n|m$.

The simplest class of supermanifolds constitute the so-called {\it split supermanifolds}.
We recall that a supermanifold $(M,\mathcal{O})$ is called {\it split} if
$\mathcal{O} = \bigwedge_{\mathcal{F}}\mathcal{G}$, where $\mathcal{G}$ is a locally free
sheaf of $\mathcal{F}$-modules on $M$. With any supermanifold $(M,\mathcal{O})$ one can
associate a split supermanifold $(M,\gr\mathcal{O})$ of the same dimension which is
called the {\it retract} of $(M,\mathcal{O})$. To construct it, let us consider the
$\mathbb{Z}_2$-graded sheaf of ideals $\mathcal{J} =
\mathcal{J}_{\bar 0}\oplus\mathcal{J}_{\bar 1}\subset\mathcal{O}$ generated by
$\mathcal{O}_{\bar 1}$. The structure sheaf of the retract is defined by
$$
\gr \mathcal{O}= \bigoplus_{p\geq 0}\gr\mathcal{O}_p,
\,\,\,\text{where} \,\,\, \gr\mathcal{O}_p=
\mathcal{J}^p/\mathcal{J}^{p+1}, \,\,\,\mathcal{J}^0:= \mathcal{O}.
$$
It can be easy shown that $\mathcal{F}\simeq
\mathcal{O}/\mathcal{J}$,  $\gr\mathcal{O}_1$ is a locally free
sheaf of $\mathcal{F}$-modules on $M$ and $\gr\mathcal{O}_p =
\bigwedge^p_{\mathcal{F}} \gr\mathcal{O}_1$. We will use the
following two locally split exact sequences:
\begin{equation}\label{sequence_stuct sheaf}
\begin{array}{c}
0\to \mathcal{J}_{\bar 0}\to \mathcal{O}_{\bar 0} \to \mathcal{F}
\to 0;\\
0\to (\mathcal{J}^2)_{\bar 1}\to \mathcal{O}_{\bar 1} \to
(\gr\mathcal{O})_1 \to 0.
\end{array}
\end{equation}
Note that a supermanifold is split iff the sequences
(\ref{sequence_stuct sheaf}) are globally split.

Let $(M,\mathcal{O})$ be a split supermanifold. Then any
$\mathbb{Z}_2$-graded locally free sheaf $\mathcal{S}=
\mathcal{S}_{\bar 0}\oplus \mathcal{S}_{\bar 1}$ of
$\mathcal{F}$-modules on $M$ gives rise to a $\mathbb{Z}_2$-graded
locally free sheaf of $\mathcal{O}$-modules $\mathcal{E}$ on
$(M,\mathcal{O})$. It is defined in the following way: $\mathcal{E}:
= \mathcal{O}\otimes_{\mathcal{F}} \mathcal{S}$. Its
$\mathbb{Z}_2$-grading is given by
\begin{equation}\label{parity of gr E}
\begin{array}{c}
\mathcal{E}_{\bar 0} = \mathcal{O}_{\bar
0}\otimes_{\mathcal{F}}\mathcal{S}_{\bar 0} +
\mathcal{O}_{\bar 1}\otimes_{\mathcal{F}} \mathcal{S}_{\bar 1},\\
\mathcal{E}_{\bar 1} =  \mathcal{O}_{\bar
0}\otimes_{\mathcal{F}}\mathcal{S}_{\bar 1} + \mathcal{O}_{\bar
0}\otimes_{\mathcal{F}}\mathcal{S}_{\bar 1}.
\end{array}
\end{equation}

Let now $\mathcal{E}=\mathcal{E}_{\bar 0}\oplus \mathcal{E}_{\bar
1}$ be a locally free sheaf of $\mathcal{O}$-modules of rang $p|q$
on an arbitrary supermanifold $(M,\mathcal{O})$. We are going to
construct a locally free sheaf of the same rank on the retract of
$(M,\mathcal{O})$. First, we note that
$\mathcal{S}:=\mathcal{E}/\mathcal{J}\mathcal{E}$ is a locally free
sheaf of $\mathcal{F}$-modules on $M$. Moreover, $\mathcal{S}$
admits the $\mathbb{Z}_2$-grading
$$
\mathcal{S} = \mathcal{S}_{\bar 0}\oplus\mathcal{S}_{\bar 1}
$$
by two locally free sheaves of $\mathcal{F}$-modules
$$
\mathcal{S}_{\bar 0}:= \mathcal{E}_{\bar 0}/(\mathcal{J}\mathcal{E})\cap\mathcal{E}_{\bar 0},
\,\, \mathcal{S}_{\bar 1}:=
\mathcal{E}_{\bar 1}/(\mathcal{J}\mathcal{E})\cap \mathcal{E}_{\bar 1}
$$
of ranks $p$ and $q$ respectively. We have the following two locally split exact sequences:
\begin{equation}\label{sequence_stuct bundle}
\begin{array}{c}
0\to \mathcal{J}\mathcal{E} \cap \mathcal{E}_{\bar 0} \to
\mathcal{E}_{(0)\bar 0}\stackrel{\alpha}{\to}
\mathcal{S}_{\bar 0} \to 0;\\
0\to \mathcal{J}\mathcal{E} \cap \mathcal{E}_{\bar 1} \to
\mathcal{E}_{(0)\bar 1} \stackrel{\beta}{\to}\mathcal{S}_{\bar 1}
\to 0,
\end{array}
\end{equation}
where $\alpha$ and $\beta$ are the natural projection maps. The
sheaf $\mathcal{E}$ possesses the filtration:
\begin{equation}\label{filtr_rassloenie}
\mathcal{E}=\mathcal{E}_{(0)}\supset \mathcal{E}_{(1)} \supset
\mathcal{E}_{(2)} \supset \ldots ,
\end{equation}
where
$$
\mathcal{E}_{(p)} = \mathcal{J}^p \mathcal{E},\,\, p\geq 1.
$$
Using this filtration, we can construct the following locally free sheaf of $\gr\mathcal{O}$-modules on the retract
$(M,\gr \mathcal{O})$:
$$
\begin{array}{c}
\gr\mathcal{E}= \bigoplus_p \gr\mathcal{E}_p,
\,\,\,\text{where} \\
 \gr\mathcal{E}_p= \mathcal{E}_{(p)} /
\mathcal{E}_{(p+1)}\simeq
\gr\mathcal{O}_p\otimes_{\mathcal{F}}\mathcal{S}.
\end{array}
$$
From $\gr\mathcal{O}= \bigwedge\gr \mathcal{O}_1$ and
$\gr\mathcal{O}_p = \bigwedge^p\gr\mathcal{O}_1$ it follows that
$$
\gr\mathcal{E} \simeq \bigwedge \gr\mathcal{O}_1
\otimes_{\mathcal{F}}\mathcal{S}.
$$
The sheaf $\gr\mathcal{E}$ we will call the {\it retract} of
$\mathcal{E}$. By definition, the sheaf $\gr\mathcal{E}$ is
$\mathbb{Z}$-graded. It possesses also the $\mathbb{Z}_{2}$-grading
given by (\ref{parity of gr E}).

Our aim now is to classify locally free sheaves of
$\mathcal{O}$-modules on a supermanifold $(M,\mathcal{O})$ which
have a fixed retract. First we formulate the well-known theorem of
Green (see [4]) which classifies complex supermanifolds
$(M,\mathcal{O}_M)$ with a given retract up to isomorphism, inducing
the identical isomorphism of reductions. The main tool used in both
classification theorems is the 1-cohomology set
$H^1(M,\mathcal{Q})$, where $\mathcal{Q}$ is a sheaf of non-abelian
groups on $M$. We denote by $\epsilon$ the unit element of
$H^1(M,\mathcal{Q})$ which corresponds to the unit 1-cocycle.

In what follows, we denote by $\mathcal{A}ut\mathcal{O}$ the sheaf
of automorphisms of the sheaf of superalgebras $\mathcal{O}$ and by
$\mathcal{A}ut^{\mathcal{R}}\mathcal{E}$ the sheaf of automorphisms
of a sheaf of $\mathcal{R}$-modules $\mathcal{E}$ on $M$, where
$\mathcal{R}$ is a sheaf of (super)algebras on $M$. The sheaf
$\mathcal{A}ut\mathcal{O}$ possesses the filtration
\begin{equation}\label{filtr_Aut_O}
\mathcal{A}ut \mathcal{O}=\mathcal{A}ut_{(0)}\mathcal{O} \supset
\mathcal{A}ut_{(2)}\mathcal{O}\supset \ldots ,
\end{equation}
where
$$
\mathcal{A}ut_{(2p)}\mathcal{O} = \{a\in\mathcal{A}ut\mathcal{O}\mid a(u)\equiv u\mod \mathcal{J}^{2p}\}.
$$
Furthermore, the group $H^0(M,\mathcal{A}ut_0\gr\mathcal{O}) \simeq
H^0(M,\mathcal{A}ut^{\mathcal{F}}\gr\mathcal{O}_1)$ acts on the
sheaf $\mathcal{A}ut \gr\mathcal{O}$ by $\Int: (a,\delta)\mapsto
a\circ \delta\circ a^{-1}$, where
$\delta\in\mathcal{A}ut\gr\mathcal{O}$ and $a\in
H^0(M,\mathcal{A}ut_0\gr\mathcal{O})$. Clearly, the group
$H^0(M,\mathcal{A}ut_0\gr\mathcal{O})$ leaves invariant the
subsheaves of groups $\mathcal{A}ut_{(2p)}\gr\mathcal{O}$. Hence
this group acts on the sets
$H^1(M,\mathcal{A}ut_{(2p)}\gr\mathcal{O})$, and the unit element
$\epsilon$ is fixed under this action.

Denote by $[(M,\mathcal{O})]$ the class of supermanifolds which are isomorphic to $(M,\mathcal{O})$.
(Here we consider complex supermanifolds up to isomorphisms
inducing the identical isomorphism of reductions.)
\medskip

\t\label{Theor_Green}[{\bf Green}] {\it Let $(M,\mathcal{O}_{\gr})$ be a split
complex supermanifold. Then
$$
\begin{array}{c}
\{[(M,\mathcal{O})] \mid \gr\mathcal{O}= \mathcal{O}_{\gr}\}
\stackrel{1:1}{\longleftrightarrow}
H^1(M,\mathcal{A}ut_{(2)}\gr\mathcal{O})/
H^0(M,\mathcal{A}ut_0\gr\mathcal{O}),
\end{array}
$$
where $(M,\mathcal{O}_{\gr})$ corresponds to $\epsilon$.
}
\medskip

\noindent {\it 2.2 Classification theorems for locally free sheaves with a
given retract}

Let $(M,\mathcal{O})$ and $(M,\mathcal{O}')$ be two supermanifolds,
$\mathcal{E}_1$ and $\mathcal{E}_2$ be locally free sheaves of
$\mathcal{O}$-modules and $\mathcal{O}'$-modules on $M$ respectively. Suppose
that $\Psi: \mathcal{O}\to\mathcal{O}'$ is a homomorphism of sheaves of
superalgebras. A homomorphism of $\mathbb{Z}_2$-graded sheaves of vector spaces
$\Phi: \mathcal{E}_1\to\mathcal{E}_2$ is called a $\Psi$-{\it morphism} if
$$
\Phi(fv)= \Psi(f)\Phi(v), \,\,f\in \mathcal{O}, \,\,v\in \mathcal{E}_1.
$$
In this case we write $\Phi = \Phi_{\Psi}$. A $\Psi$-{\it morphism}
$\Phi: \mathcal{E}\to\mathcal{E}$ is called a $\Psi$-{\it
isomorphism} if $\Phi$ is invertible. A $\Psi$-isomorphism $\Phi:
\mathcal{E}\to\mathcal{E}$ we also will call a $\Psi$-{\it
automorphism} of $\mathcal{E}$. A homomorphism (isomorphism) of
$\mathbb{Z}_2$-graded sheaves of vector spaces $\Phi:
\mathcal{E}_1\to\mathcal{E}_2$ will be called a {\it quasi-morphism}
({\it quasi-isomorphism}) if it is a $\Psi$-{\it morphism}
($\Psi$-{\it isomorphism}) for a certain $\Psi$. The sheaves
$\mathcal{E}_1$ and $\mathcal{E}_2$ will be called {\it
quasi-isomorphic} if it exists a quasi-isomorphism $\Phi:
\mathcal{E}_1\to\mathcal{E}_2$. A quasi-isomorphism
$\mathcal{E}\to\mathcal{E}$ will be called a {\it
quasi-automorphism} of $\mathcal E$. We will study the sheaf
$\mathcal{QA}ut\mathcal{E}$, where
\begin{equation}
\begin{split}
\mathcal{QA}ut\mathcal{E}(U) =\{\Phi\,\mid \,\Phi\,\,\,\text{is a
quasi-automorphism of}\,\,\mathcal{E}|_U\}
\end{split}
\end{equation}
for each open subset $U\subset M$. One verifies easily that $\Phi_{\Psi}\circ
\Theta_{\Upsilon}$, where $\Phi_{\Psi},\, \Theta_{\Upsilon}\in
\mathcal{QA}ut\mathcal{E}$, is a $\Psi\circ \Upsilon$-morphism. It follows that
$\mathcal{QA}ut\mathcal{E}$ is a sheaf of groups. It possesses the double
filtration by the subsheaves
$$
\begin{array}{c}
\mathcal{QA}ut_{(p)(q)} \mathcal{E}:= \{\Phi_{\Psi}\in \mathcal{QA}ut\mathcal{E}\mid
\Phi_{\Psi}(v)\equiv v \mod \mathcal{E}_{(p)},\,
\Psi(f)\equiv f\mod \mathcal{J}^q\\
\text{for}\,\, v\in \mathcal{E}, f\in \mathcal{O} \},\,\,p,q\geq 0.
\end{array}
$$
We also define the following subsheaves:
\begin{equation}
\begin{split}
\mathcal{QA}ut_0(\gr\mathcal{E}):= \{\Phi_{\Psi} \mid \Phi_{\Psi}\in \mathcal{QA}ut (\gr\mathcal{E}), \,\,
\Phi_{\Psi}\, \text{preserves the $\mathbb{Z}$-grading of}\, \gr\mathcal{E}\}.
\end{split}
\end{equation}

\begin{equation}
\mathcal{A}ut^{\mathcal{F}}_{\bar 0}\mathcal{S}:= \{\Phi \mid
\Phi\in\mathcal{A}ut^{\mathcal{F}}\mathcal{S}, \,\, \Phi\,
\text{preserves the $\mathbb{Z}_2$-grading of}\, \mathcal{S}\},
\end{equation}
where $\mathcal{S}$ is a $\mathbb{Z}_2$-graded sheaf of
$\mathcal{F}$-modules.

\medskip

\lem\label{Lemma Aut gr E} {\it We have an isomorphism of sheaves of groups
$$
\mathcal{QA}ut_0(\gr\mathcal{E}) \simeq {\mathcal
A}ut^{\mathcal{F}}(\gr\mathcal{O}_1)
\times\mathcal{A}ut^{\mathcal{F}}_{\bar 0}\mathcal{E}_{\red}.
$$
}
\medskip

\noindent{Proof.} Let us define the mapping
$$
\Theta: {\mathcal A}ut^{\mathcal{F}}(\gr\mathcal{O}_1)
\times\mathcal{A}ut^{\mathcal{F}}_{\bar 0}\mathcal{E}_{\red} \to
\mathcal{QA}ut_0(\gr\mathcal{E})
$$
by
$$
(\psi,\Phi)\mapsto \Phi_{\wedge\psi}, \,\, \psi\in {\mathcal
A}ut^{\mathcal{F}}(\gr\mathcal{O}_1), \,\, \Phi\in
\mathcal{A}ut^{\mathcal{F}}_{\bar 0}\mathcal{E}_{\red},
$$
where
$$
\Phi_{\wedge\psi}(hv):= \wedge\psi(h)
\Phi(v)
$$
for $h\in\gr\mathcal{O}$, $v\in\mathcal{E}_{\red}$ and $\wedge\psi$ is the automorphism of the sheaf
$\gr\mathcal{O}$ induced by $\psi$. This is a homomorphism of sheaves of groups. In fact, suppose that another pair
$(\psi',\Phi')$, where $\psi'\in
{\mathcal
A}ut^{\mathcal{F}}(\gr\mathcal{O}_1)$, $\Phi'\in
\mathcal{A}ut^{\mathcal{F}}_{\bar 0}\mathcal{E}_{\red}$, is given. Then we have
$$
\begin{array}{c}
(\Phi_{\wedge\psi}\circ\Phi'_{\wedge\psi'})(hv)= \Phi_{\wedge\psi}(\wedge\psi'(h)\Phi'_{\wedge\psi'}(v)) =
\wedge\psi(\wedge\psi'(h))\Phi_{\wedge\psi}(\Phi'_{\wedge\psi'}(v)) = \\
(\Phi\circ\Phi'_{\wedge\psi\circ \wedge\psi'})(hv)
\end{array}
$$
for $h\in\gr\mathcal{O},\, v\in\mathcal{E}_{\red}$.

Let us prove that $\Ker\Theta = (\id,\id)$. Suppose that $\Theta(\psi,\Phi) =
\id$. Then $\Phi_{\wedge\psi}(hv) = \wedge\psi(h)\Phi(v) = hv$ for all $h\in\gr\mathcal{O},\,
v\in\mathcal{E}_{\red}$. Putting $h = 1$, we see that $\Phi(v) = v$, i.e.,
$\Phi = \id$. Since $\mathcal{E}_{\red}$ is locally free, this implies that
$\wedge\psi(h) = h$, therefore, $\psi = \id$. Thus, the homomorphism $\Theta$ is injective.

Let us now prove that it is surjective. Let $\Phi_{\Psi}\in
\mathcal{QA}ut_0(\gr \mathcal{E})$ be given. Let us show that
$\Phi_{\Psi}\in \Im \Theta$. Since
$\Phi_{\Psi}|_{\mathcal{E}_{\red}}: \mathcal{E}_{\red}\to
\mathcal{E}_{\red}$ and $\Phi_{\Psi}$ preserves the
$\mathbb{Z}_{2}$-grading of $\gr \mathcal{E}$, we have $\Phi:=
\Phi_{\Psi}|_{\mathcal{E}_{\red}}\in
\mathcal{A}ut^{\mathcal{F}}_{\bar 0}\mathcal{E}_{\red}$.
Furthermore, if $h\in\gr\mathcal{O}_p$ and $v\in
\mathcal{E}_{\red}$, then
$$
\Phi_{\Psi}(hv) = \Psi(h)\Phi(v)\in \gr\mathcal{E}_p.
$$
It follows that $\Psi(h)\in \gr\mathcal{O}_p$, and hence $\Psi$
preserves the $\mathbb{Z}$-grading of $\gr\mathcal{O}$. We have
$\psi = \Psi|_{\gr\mathcal{O}_1}\in {\mathcal
A}ut^{\mathcal{F}}(\gr\mathcal{O}_1)$ and $\wedge\psi = \Psi$. The
proof is complete.$\Box$

\medskip

We will use the above notation, fixing a split complex supermanifold
$(M,\mathcal{O}_{\gr})$ and a $\mathbb{Z}_2$-graded locally free
sheaf of $\mathcal{F}$-modules $\mathcal{S}$ on $M$. Our aim is to
classify locally free sheaves $\mathcal{E}$ of $\mathcal{O}$-modules
on complex supermanifolds $(M,\mathcal{O})$ with retract
$(M,\mathcal{O}_{\gr})$, whose retract $\gr\mathcal{E}$ coincides
with $\mathcal{E}_{\gr} =
\mathcal{O}_{\gr}\otimes_{\mathcal{F}}\mathcal{S}$.

The group $H^0(M,\mathcal{QA}ut_0 \mathcal{E}_{\gr})$ acts on the sheaf $\mathcal{QA}ut\mathcal{E}_{\gr}$ by
the automorphisms $\delta\mapsto a\circ \delta \circ a^{-1}$, where $a\in
H^0(M,\mathcal{QA}ut_0 \mathcal{E}_{\gr})$ and $\delta\in\mathcal{QA}ut \mathcal{E}_{\gr}$.
It is easy to see that this action leaves invariant the subsheaves
$\mathcal{QA}ut_{(p)(q)}\mathcal{E}_{\gr}$ and hence induces an action of
$H^0(M,\mathcal{QA}ut_0 \mathcal{E}_{\gr})$ on the cohomology set
$H^1(M,\mathcal{QA}ut_{(p)(q)}\mathcal{E}_{\gr})$.

If $\phi:M\to N$ is a holomorphic map of manifolds and $p:\mathbb{E}\to N$ is a
vector bundle, we may define the pullback bundle $\phi^*(\mathbb{E})$ on $N$.
The corresponding to $\phi^*(\mathbb{E})$ sheaf is
$\mathcal{F}_M\otimes_{\phi^*(\mathcal{F}_N)} \phi^*(\mathcal{E})$, where
$\mathcal{E}$ is the sheaf of sections corresponding to $\mathbb{E}$,
$\mathcal{F}_M$ and $\mathcal{F}_N$ are the sheaves of holomorphic functions on
$M$ and $N$ respectively.
Let $\pi:(M,\mathcal{O}_M) \to (N,\mathcal{O}_N)$ be a morphism of two
supermanifolds and $\mathcal{E}$ be a locally free sheaf of
$\mathcal{O}_N$-modules on $N$ of rang $p|q$. Similarly, we can define the
sheaf $\mathcal{O}_M\otimes_{\pi_{\red}^*(\mathcal{O}_N)}
\pi_{\red}^*(\mathcal{E})$. This sheaf is a locally free sheaf of
$\mathcal{O}_M$-modules on $M$ of rang $p|q$, since
$$
\mathcal{O}_M\otimes_{\pi_{\red}^*(\mathcal{O}_N)}
\pi_{\red}^*(\mathcal{O}_N)\simeq \mathcal{O}_M.
$$
Sometimes we will denote the sheaf
$\mathcal{O}_M\otimes_{\pi_{\red}^*(\mathcal{O}_N)}
\pi_{\red}^*(\mathcal{E})$ by $\widetilde{\pi}(\mathcal{E})$.

Let us consider the special case $(M,\mathcal{O}_M) =
(N,\mathcal{O}_N)$, $\pi=(\id,\pi^*)$ and $\pi^*\in
H^0(M,\mathcal{A}ut \mathcal{O}_M)$. We have
$$
\widetilde{\pi}(\mathcal{E})=
\mathcal{O}_M\otimes_{\id^*(\mathcal{O}_N)} \id^*(\mathcal{E}) =
\mathcal{O}_M\otimes_{\mathcal{O}_N} \mathcal{E}.
$$
The sheaves $\widetilde{\pi}(\mathcal{E})$ and $\mathcal{E}$ are
$(\pi^*)^{-1}$-isomorphic, the $(\pi^*)^{-1}$-isomorphism is given
by $f\otimes s\mapsto (\pi^*)^{-1}(f)s$, where $f\in \mathcal{O}_M$
and $s\in \mathcal{E}$. Let $\Phi_{\Psi^*}:\mathcal{E}\to
\mathcal{E}'$ be an $\Psi^*$-isomorphism of two locally free sheaves
of $\mathcal{O}_M$-modules on $M$. We put $\Psi:=(\id, \Psi^*)$. We
see that
 $\widetilde{\Psi}(\mathcal{E})$ and $\mathcal{E}'$ are $\id$-isomorphic.

Furthermore, let us consider the sheaf $\mathcal{A}ut^{\mathcal{O}}\mathcal{E}$ of
automorphisms of the $\mathcal{O}$-modules sheaf $\mathcal{E}$. It possesses the filtration:
$$
\mathcal{A}ut^{\mathcal{O}}\mathcal{E}=\mathcal{A}ut^{\mathcal{O}}_{(0)}\mathcal{E}
\supset\mathcal{A}ut^{\mathcal{O}}_{(1)}\mathcal{E}\supset
\ldots,
$$
where
$$
\mathcal{A}ut_{(p)}^{\mathcal{O}}\mathcal{E}:=\{a\in \mathcal{A}ut^{\mathcal{O}}\mathcal{E} \mid a(v)\equiv v\mod \mathcal{E}_{(p)}\}, \,\,p\geq 0.
$$
The group $H^0(M,\mathcal{A}ut_{0}^{\mathcal{O}}\gr\mathcal{E})\simeq H^0(M,\mathcal{A}ut_{\bar 0}^{\mathcal{F}}\mathcal{E}_{\red})$  acts on the sheaf
$\mathcal{A}ut^{\mathcal{O}}\gr\mathcal{E}$ by $\delta\mapsto a\circ \delta \circ a^{-1}$,
where $a\in H^0(M,\mathcal{A}ut_{0}^{\mathcal{O}}\gr\mathcal{E})$ and $\delta\in
\mathcal{A}ut^{\mathcal{O}}\gr\mathcal{E}$. It is easy to see that this action leaves
 the subsheaves $\mathcal{A}ut_{(p)}^{\mathcal{O}}\gr\mathcal{E}$ invariant and hence
induces an action of $H^0(M,\mathcal{A}ut_{0}^{\mathcal{O}}\gr\mathcal{E})$ on the cohomology set
$H^1(M,\mathcal{A}ut_{(p)}^{\mathcal{O}}\gr\mathcal{E})$.

We have the exact sequence of sheaves of groups
$$
\id\to \mathcal{A}ut^{\mathcal{O}}\mathcal{E} \to \mathcal{QA}ut\mathcal{E}
\to \mathcal{A}ut\mathcal{O}\to\id,
$$
where the first homomorphism is the natural embedding (an automorphism of
$\mathcal{A}ut^{\mathcal{O}}\mathcal{E}$ is regarded as an $\id$-morphism) and the second
one, say $F:\mathcal{QA}ut \mathcal{E} \to \mathcal{A}ut \mathcal{O}$, is
defined by $\Phi_{\Psi} \mapsto \Psi$.
Note that $F(\mathcal{QA}ut_{(p)(q)} \mathcal{E})\subset
\mathcal{A}ut_{(q)}\mathcal{O}$ and in the case $\mathcal{E}=\gr\mathcal{E}$
the restriction $F|\mathcal{QA}ut_0\gr\mathcal{E}$ coincides with the
natural projection
$$
\mathcal{QA}ut_0(\mathcal{E}_{\gr})\simeq \mathcal{A}ut_0\gr\mathcal{O}\times \mathcal{A}ut^{\mathcal{F}}_{\bar 0}(\mathcal{E}_{\red})
\to\mathcal{A}ut_0\gr\mathcal{O}
$$
(see Lemma 1).

The homomorphism $F$ commutes with the actions of $H^0(M,\mathcal{QA}ut_{0}\gr\mathcal{E})$
and $H^0(M,\mathcal{A}ut_{0}\gr\mathcal{O})$ on $\mathcal{QA}ut_{(p)(q)}(\gr\mathcal{E})$ and
$\mathcal{A}ut_{(q)}(\gr\mathcal{O})$ respectively. More precisely,
$$
F(a\circ \delta \circ a^{-1}) = F(a)\circ F(\delta) \circ F(a^{-1}),
$$
where $a\in H^0(M,\mathcal{QA}ut_{0}\gr\mathcal{E})$ and
$\delta\in\mathcal{QA}ut\gr\mathcal{E}$. It follows that $F$ induces the map
of sets
$$
\begin{array}{ll}
\widetilde{F}: H^1(M,\mathcal{QA}ut_{(1)(2)}
\gr\mathcal{E})/&\!\!\!\!\!H^0(M,\mathcal{QA}ut_{0}\gr\mathcal{E})\to \\
& H^1(M,\mathcal{A}ut_{(2)}\gr\mathcal{O})/H^0(M,\mathcal{A}ut_{0}\gr\mathcal{O}).
\end{array}
$$

Let $\Phi_{\Psi}:\mathcal{E}_1\to \mathcal{E}_2$ be a $\Psi$-morphism of
locally free sheaves of $\mathcal{O}$-modules. Since
$\Psi(\mathcal{J}^p)\subset \mathcal{J}^p$, we see that
$\Phi_{\Psi}((\mathcal{E}_1)_{(p)})\subset (\mathcal{E}_2)_{(p)}$, $p\geq 0$.
We denote by $\gr(\Phi_{\Psi}): \gr\mathcal{E}_1\to\gr \mathcal{E}_2$ the
induced morphism. Let $\mathcal{E}$ be a locally free sheaf of $\mathcal{O}$-modules on $M$.
Denote
$$
[\mathcal{E}] = \{ \mathcal{E}' \mid \mathcal{E}' \,\,\, \text{is quasi-isomorphic to} \,\,\, \mathcal{E}  \}.
$$

\medskip

\t\label{Aut_E} {\it Let $(M,\mathcal{O}_{\gr})$ be a split supermanifold,
$\mathcal{S} = \mathcal{S}_{\bar 0}\oplus \mathcal{S}_{\bar 1}$ be a
$\mathbb{Z}_2$-graded locally free sheaf of $\mathcal{F}$-modules on $M$ and
$\mathcal{E}_{\gr} = \mathcal{O}_{\gr}\otimes_{\mathcal{F}}\mathcal{S}$.

1)We have a bijection
$$
\begin{array}{c}
\{[\mathcal{E}] \mid \gr\mathcal{O}= \mathcal{O}_{\gr}, \,\, \gr\mathcal{E} =
\mathcal{E}_{\gr}\}\stackrel{1:1}{\longleftrightarrow}
 H^1(M,\mathcal{QA}ut_{(1)(2)}\mathcal{E}_{\gr})/H^0(M,\mathcal{QA}ut_{0}\mathcal{E}_{\gr}).
\end{array}
$$
The unit $\epsilon\in H^1(M,\mathcal{QA}ut_{(1)(2)}\mathcal{E}_{\gr})$ is fixed
with respect to the action of the group
$H^0(M,\mathcal{QA}ut_{0}\mathcal{E}_{\gr})$.

2) Let $a\in H^1(M,\mathcal{A}ut_{(2)}\mathcal{O}_{\gr})/H^0(M,\mathcal{A}ut_{0}\mathcal{O}_{\gr})$.
Then there is a bijection between elements of the set
$\widetilde{F}^{-1}(a)$ and classes of isomorphic locally free sheaves on
supermanifolds which are
contained in $[(M,\mathcal{O})]$.

 }

\medskip

\noindent{Proof.}  Let $\mathcal{E}$ be a locally free sheaf of
$\mathcal{O}$-modules on $(M,\mathcal{O})$ and $\mathcal{U}=\{ U_i\}$ be an
open covering of $M$ such that (\ref{sequence_stuct sheaf}) and
(\ref{sequence_stuct bundle}) are split over $U_i$ and $\mathcal{E}|_{U_i}$ are
free. In this case $(\gr\mathcal{E})|_{U_i}$ are free sheaves of
$(\gr\mathcal{O})$-modules, too. We fix local bases $(\hat{e}^i_j)$ and
$(\hat{f}^i_k)$ of the sheaves of $\mathcal{F}$-modules
$(\mathcal{E}_{\red})_{\bar 0}|_{U_i}$ and $(\mathcal{E}_{\red})_{\bar
1}|_{U_i}$, $U_i\in \mathcal{U}$, respectively.

 We are going to define an isomorphism $\delta_i:\mathcal{E}|_{U_i}\to(\gr\mathcal{E})|_{U_i}$. Let
$e^i_j\in \mathcal{E}_{(0)\bar 0}$ such that $\alpha(e^i_j)=\hat{e}^i_j$ and
$f^i_k\in \mathcal{E}_{(0)\bar 1}$ such that $\beta(f^i_k)=\hat{f}^i_k$. Then
$(e^i_j,f^i_k)$ is a local basis of $\mathcal{E}|_{U_i}$. A splitting of
(\ref{sequence_stuct sheaf}) determines local isomorphisms
$\sigma_i:\mathcal{O}|_{U_i} \to \gr\mathcal{O}|_{U_i}$. We put
$$
\delta_i(\sum h_je^i_j+ \sum g_kf^i_k)=  \sum
\sigma_i(h_j)\hat{e}^i_j+ \sum\sigma_i(g_k)\hat{f}^i_k,
\,\,\,h_j,g_k\in \mathcal{O}.
$$
Obviously, $\delta_i$ is an isomorphism. We put $\gamma_{ij}:=
\sigma_i\circ\sigma_j^{-1}$ and $(g_{ij})_{\gamma_{ij}}:= \delta_i\circ
\delta_j^{-1}$. It is clear that $(\gamma_{ij})\in Z^1(\mathcal{U},
\mathcal{A}ut_{(2)}(\gr\mathcal{O}))$ and
$$
((g_{ij})_{\gamma_{ij}})\in Z^1(\mathcal{U}, \mathcal{QA}ut_{(1)(2)}(\gr\mathcal{E})).
$$

Conversely, if $((g_{ij})_{\gamma_{ij}})\in
Z^1(\mathcal{U},\mathcal{QA}ut_{(1)(2)}(\gr\mathcal{E}))$, we can construct a
locally free sheaf of $\mathcal{O}$-modules on $(M,\mathcal{O}(\gamma_{ij}))$,
where $(M,\mathcal{O}(\gamma_{ij}))$ is the supermanifold corresponding to the
cocycle $(\gamma_{ij})\in Z^1(\mathcal{U},\mathcal{A}ut_{(2)}\gr \mathcal{O})$
by the Green Theorem. Indeed, we have to identify $\gr\mathcal{E}|_{U_i}$ with
$\gr\mathcal{E}|_{U_j}$ over $U_i\cap U_j$ using $(g_{ij})_{\gamma_{ij}}$.

The standard calculation shows that if two cocycles $((g_{ij})_{\gamma_{ij}}) $
and $((g_{ij}')_{\gamma_{ij}'})$ are cohomological, then the corresponding
locally free sheaves of $\mathcal{O}$-modules are quasi-isomorphic and this
quasi-isomorphism denoted by $\Phi_{\Psi}$ has the property
$\gr(\Phi_{\Psi})=\id_{\id}$. Conversely, if $\Phi_{\Psi}:\mathcal{E}\to
\mathcal{E}'$ is a quasi-isomorphism of locally free sheaves of
$\mathcal{O}$-modules such that $\gr(\Phi_{\Psi})=\id_{\id}$, then the
corresponding cocycles are cohomological.

Let $\mathcal{E}$ and $\mathcal{E}'$ be two locally free sheaves of
$\mathcal{O}$-modules on $(M,\mathcal{O})$ such that $\gr
\mathcal{E}=\gr \mathcal{E}'=\mathcal{E}_{\gr}$. Assume that
$\Phi_{\Psi}:\mathcal{E}\to \mathcal{E}'$ is an isomorphism. Then
$\gr(\Phi_{\Psi})\in H^0(M,\mathcal{QA}ut_0 \gr \mathcal{E})$. Suppose that
$\mathcal{E}$ corresponds to $(g_{ij})_{\gamma_{ij}}=\delta_i\circ
\delta_j^{-1}$, where $\gamma_{ij}=\sigma_i\circ \sigma_j^{-1}$, and
$\mathcal{E}'$ corresponds to
$(g_{ij}')_{\gamma_{ij}'}=\delta'_i\circ (\delta'_j)^{-1}$, where
$\gamma'_{ij}=\sigma'_i\circ (\sigma'_j)^{-1}$. There exist
isomorphisms
$(\widetilde{\Phi}_i)_{\widetilde{\Psi}_i}:\gr\mathcal{E}|_{U_i}\to
\gr\mathcal{E}|_{U_i}$ such that the following diagram is
commutative:
$$
\begin{CD}
\gr\mathcal{E}|_{U_i} &
@>{(\widetilde{\Phi}_i)_{\widetilde{\Psi}_i}}>>&
\gr\mathcal{E}|_{U_i}\\
@A{\delta_i}AA && @AA{\delta_i'}A\\
\mathcal{E}|_{U_i} & @>{\Phi_{\Psi}}>>& \mathcal{E}|_{U_i}
\end{CD}.
$$
Since $\gr \delta_i = \gr \delta'_i$, it follows that $\gr
((\widetilde{\Phi}_i)_{\widetilde{\Psi}_i}) = \gr (\Phi_{\Psi})$ and
hence
$$
(\Theta_i)_{\Omega_i}:=\gr(\Phi_{\Psi})^{-1} \circ
(\widetilde{\Phi}_i)_{\widetilde{\Psi}_i}\in \mathcal{QA}ut_{(1)(2)}
\gr\mathcal{E}.
$$
Further, we have
$$
\begin{array}{c}
(g_{ij}')_{\gamma_{ij}'} = \delta'_i \circ (\delta'_j)^{-1} =
(\widetilde{\Phi}_i)_{\widetilde{\Psi}_i} \circ \delta_i \circ
(\Phi_{\Psi})^{-1} \circ \Phi_{\Psi} \circ \delta_j^{-1} \circ
((\widetilde{\Phi}_j)_{\widetilde{\Psi}_j})^{-1}=\\
(\widetilde{\Phi}_i)_{\widetilde{\Psi}_i} \circ
(g_{ij})_{\gamma_{ij}} \circ
((\widetilde{\Phi}_j)_{\widetilde{\Psi}_j})^{-1}=
\gr(\Phi_{\Psi})\circ (\Theta_i)_{\Omega_i}\circ
(g_{ij})_{\gamma_{ij}} \circ (\Theta_j^{-1})_{\Omega_j^{-1}}\circ
\gr(\Phi_{\Psi})^{-1}.
\end{array}
$$
Hence, the cohomology classes corresponding to
$(g_{ij})_{\gamma_{ij}}$ and $(g_{ij}')_{\gamma_{ij}'}$ belong to
the same orbit of the group $H^0(M,\mathcal{QA}ut_0
\mathcal{E}_{\gr})$.

Conversely, assume that $b\in H^0(M,\mathcal{QA}ut_0
\mathcal{E}_{\gr})$ and $(g_{ij}')_{\gamma_{ij}'}=b \circ
(g_{ij})_{\gamma_{ij}} \circ b^{-1}$. Then $\delta'_i \circ
(\delta'_j)^{-1} = b \circ \delta_i \circ \delta_j^{-1} \circ
b^{-1}$ and we can define the isomorphism $\Gamma: \mathcal{E}\to
\mathcal{E}'$ by $\Gamma|_{U_i}:= (\delta'_i)^{-1} \circ b\circ
\delta_i$, where $\mathcal{E}$ and $\mathcal{E}'$ correspond to
$(g_{ij})_{\gamma_{ij}}$ and $(g_{ij}')_{\gamma_{ij}'}$
respectively.

Let $a\in H^1(M,\mathcal{A}ut_{(2)} \mathcal{O}_{\gr})/H^0(M,\mathcal{A}ut_{0} \mathcal{O}_{\gr})$.
By Theorem \ref{Theor_Green}  we may assign to
each $a$ the class of isomorphic supermanifolds $[(M,\mathcal{O})]$.
From the proof of Theorem \ref{Aut_E} it follows that there is a
bijection between elements of the set $\widetilde{F}^{-1}(a)$ and
classes of isomorphic locally free sheaves on  supermanifolds which
are contained in $[(M,\mathcal{O})]$.$\Box$

\medskip

\noindent {\it 2.3 A classification theorem for locally free sheaves on a split
supermanifold}

Denote by $[\mathcal{E}]_{\id}$ the class of $\id$-isomorphic (i.e.,
isomorphic) to $\mathcal{E}$ locally free sheaves of $\mathcal{O}$-modules on a
split complex supermanifold $(M,\mathcal{O})$.

\medskip

\t\label{Aut_E_(gr)} {\it Let $(M,\mathcal{O})$ be a split supermanifold,
$\mathcal{S} = \mathcal{S}_{\bar 0}\oplus\mathcal{S}_{\bar 1}$ be a
$\mathbb{Z}_2$-graded locally free sheaf of $\mathcal{F}$-modules on $M$ and
$\mathcal{E}_{\gr} = \mathcal{O}\otimes_{\mathcal{F}} \mathcal{E}_{\red}$. Then
$$
\begin{array}{c}
\{[\mathcal{E}]_{\id} \mid \gr \mathcal{E}= \mathcal{E}_{\gr}\}
\stackrel{1:1}{\longleftrightarrow} H^1(M,\mathcal{A}ut^{\mathcal{O}}_{(1)}\mathcal{E}_{\gr})/H^0(M,\mathcal{A}ut^{\mathcal{O}}_{ 0}\mathcal{E}_{\gr}).
\end{array}
$$
Moreover, the unit $\epsilon\in H^1(M,\mathcal{A}ut_{(1)}^{\mathcal{O}}\mathcal{E}_{\gr})$ is
a fixed point with respect to the action of the group $H^0(M,\mathcal{A}ut^{\mathcal{O}}_{ 0}\mathcal{E}_{\gr})$.
}

\medskip

\noindent{Proof.} Let us use the notations from the proof of Theorem 2. Since
$(M,\mathcal{O})$ is split, we may assume that $\sigma_i=\sigma|_{U_i}$, where $\sigma$ is determined by a global splitting of (\ref{sequence_stuct sheaf}). It follows
that the cocycle $(g_{ij})$ lies in $Z^1(\mathcal{U},\mathcal{A}ut_{(1)}^{\mathcal{O}}
\mathcal{E}_{\gr})$. The further proof is similar to the proof the Theorem 2.
$\Box$

\bigskip

\begin{center}
    {\bf 3. Locally free sheaves of modules on projective superspaces}
\end{center}

In this subsection we will discuss two remarkable theorems about
locally free sheaves on projective spaces, proved by Barth -- Van de
Ven -- Tyurin and Birkhoff -- Grothendieck, in the super-context.

\smallskip

\noindent {\it 3.1 Exact sequences corresponding to} $\mathcal{A}ut^{\mathcal{O}}\mathcal{E}$

\smallskip

Let $(M,\mathcal{O})$ be a split complex supermanifold and
$\mathcal{E}$ be a locally free sheaf of $\mathcal{O}$-modules on
$M$. Denote by $\mathcal{E}nd^{\mathcal{O}}\mathcal{E}$ the sheaf of
$\mathcal{O}$-endomorphisms of $\mathcal{E}$. This sheaf possesses
the filtration
$$
\mathcal{E}nd^{\mathcal{O}} \mathcal{E} =
\mathcal{E}nd^{\mathcal{O}}_{(0)}\mathcal{E} \supset
\mathcal{E}nd^{\mathcal{O}}_{(1)} \mathcal{E} \supset \ldots ,
$$
$$
\mathcal{E}nd^{\mathcal{O}}_{(p)} \mathcal{E}:= \{ A\in
\mathcal{E}nd^{\mathcal{O}} \mathcal{E} \mid
A(\mathcal{E}_{(q)})\subset \mathcal{E}_{(q+p)}\,\,\text{for
all}\,\,q\geq 0\}.
$$
The map
$$
\exp: \mathcal{E}nd^{\mathcal{O}}_{(p)} \mathcal{E} \to
\mathcal{A}ut^{\mathcal{O}}_{(p)} \mathcal{E},
$$
given by the usual $\exp$-series is a bijection of sheaves of sets
for all $p\geq 1$ due to the fact that $\log = (\exp)^{-1}$ is well
defined. In general it is not a homomorphism of sheaves of groups.
We may define the map
$$
\lambda_p:\mathcal{A}ut^{\mathcal{O}}_{(p)} \mathcal{E}\to
\mathcal{E}nd^{\mathcal{O}}_{(p)}\mathcal{E}/\mathcal{E}nd^{\mathcal{O}}_{(p+1)}
\mathcal{E}, \,\,p\geq 1,
$$
 given by
$$
a\mapsto A+\mathcal{E}nd^{\mathcal{O}}_{(p+1)},
\,\,\text{where}\,\,a=\exp(A).
$$
This map is surjective and $\Ker
\lambda_p=\mathcal{A}ut^{\mathcal{O}}_{(p+1)}\mathcal{E}$. Clearly,
it is a homomorphism of sheaves of groups. We will also consider the
subsheaves of $\mathcal{E}nd^{\mathcal{O}}\gr\mathcal{E}$
$$
\mathcal{E}nd^{\mathcal{O}}_{p} \gr\mathcal{E}:= \{ A\in
\mathcal{E}nd^{\mathcal{O}} \gr\mathcal{E} \mid
A(\gr\mathcal{E}_{q})\subset \gr\mathcal{E}_{p+q}\}, \,\,p\geq 0.
$$
Then
$$
\mathcal{E}nd^{\mathcal{O}}_{(p)} \gr\mathcal{E}=\bigoplus_{q\geq p}
\mathcal{E}nd^{\mathcal{O}}_{q} \gr\mathcal{E}.
$$
It follows that
$$
\mathcal{E}nd^{\mathcal{O}}_{(p)}\gr\mathcal{E}/\mathcal{E}nd^{\mathcal{O}}_{(p+1)}
\gr\mathcal{E} \simeq \mathcal{E}nd^{\mathcal{O}}_{p}
\gr\mathcal{E}.
$$
Hence, we get the exact sequence
\begin{equation}\label{sequense Aut split}
0\to \mathcal{A}ut^{\mathcal{O}}_{(p+1)}\gr \mathcal{E} \to
\mathcal{A}ut^{\mathcal{O}}_{(p)}
\gr\mathcal{E}\stackrel{\lambda_p}{\to}
\mathcal{E}nd^{\mathcal{O}}_{p} \gr\mathcal{E}\to 0, \,\,p\geq 1.
\end{equation}

The following lemma gives a description of the sheaf
$\mathcal{E}nd_{p}\gr\mathcal{E},\, p\geq 1$, in terms of the
sheaves $\mathcal{O}$ and $\mathcal{E}_{\red}$.

\medskip

\lem\label{Aut_E_(p)=} {\it We have
$$
\begin{array}{rl}
\mathcal{E}nd^{\mathcal{O}}_{p} \gr\mathcal{E}  \simeq&\left\{
                                         \begin{array}{ll}
                                           \mathcal{O}_{p} \otimes
((\mathcal{E}_{\red})_{\bar 0} \otimes (\mathcal{E}_{\red})_{\bar
1}^*  \oplus (\mathcal{E}_{\red})_{\bar 1} \otimes
(\mathcal{E}_{\red})_{\bar 0}^*), & \hbox{$p$ is odd;} \\
                                        \mathcal{O}_{p} \otimes
((\mathcal{E}_{\red})_{\bar 0} \otimes (\mathcal{E}_{\red})_{\bar
0}^*  \oplus (\mathcal{E}_{\red})_{\bar 1} \otimes
(\mathcal{E}_{\red})_{\bar 1}^*)   , & \hbox{$p$ is even.}
                                         \end{array}
                                       \right.
\end{array}
$$
}

\noindent{Proof.} Firstly, note that an endomorphism
$A\in\mathcal{E}nd_{p}(\gr\mathcal{E})$ is determined by its
restriction $A|_{\gr\mathcal{E}_{0}}$. Secondly, $A|_{\gr
\mathcal{E}_{0}} : \gr \mathcal{E}_{0} \to \gr \mathcal{E}_{p}$ is
an $\mathcal{F}$-linear map preserving parity (\ref{parity of gr
E}). The result follows from the relation $\gr\mathcal{E}_{q}\simeq
\gr\mathcal{O}_{q}\otimes\mathcal{E}_{\red}$. $\Box$

\smallskip

Now we can recover the following well-known result, see \cite{Man,Rempel}:

\smallskip
\l\label{smooth loc free sheaves}{\it Let $(M,\mathcal{O})$ be a smooth
supermanifold and $\mathcal{E}$ be a locally free sheaf of
$\mathcal{O}$-modules on $M$. Then $\mathcal{E}\simeq
\mathcal{O}\otimes_{\mathcal{F}} \mathcal{E}_{\red}$. }

\medskip

\noindent{\it Proof.} Indeed, $(M,\mathcal{O})$ is split by the
Batchelor Theorem. In this case
$$
H^1(M,\mathcal{E}nd_{p}^{\mathcal{O}}\gr\mathcal{E})=\{0\}
$$
by Lemma \ref{Aut_E_(p)=}. Hence
$$
H^1(M,\mathcal{A}ut_{(1)}^{\mathcal{O}}\gr\mathcal{E})=\{\epsilon\},
$$
and our assertion follows from the Theorem \ref{Aut_E_(gr)}.$\Box$

\bigskip

\noindent {\it 3.2 The Barth -- Van de Ven -- Tyurin Theorem for
supermanifolds}

\smallskip
Let us briefly recall the classical Barth -- Van de Ven -- Tyurin
Theorem. Consider the sequence of complex projective spaces
$$
\mathbb{CP}^{1}\stackrel{\varphi_{1}}{\longrightarrow}
\mathbb{CP}^{2}\stackrel{\varphi_{2}}{\longrightarrow} \ldots,
$$
where $\varphi_{i}$ are standard embeddings. (The inductive limit of this
sequence is also called the {\it complex projective ind-space}
$\mathbb{CP}^{\infty}$ (see \cite{Don_Penkov, Tyurin} and more detailed
\cite{Kumar}.) We consider collections $E=\{E_N\}_{N\geq 1}$ of holomorphic
vector bundles $E_N$ of a finite rank over $\mathbb{CP}^{N}$, $N\geq 1$, such
that $\widetilde{\varphi}_{N} (E_{N+1}) =E_N$. (Such collections are also
called {\it vector bundles over $\mathbb{CP}^{\infty}$}.) If $E=\{E_N\}_{N\geq
1}$ and $E'=\{E'_N\}_{N\geq 1}$ are two such collections, then the collection
$E\oplus E':=\{E_N\oplus E'_N\}_{N\geq 1}$ is called the {\it direct sum} of
$E$ and $E'$. A {\it morphism of collections} $f:E\to E'$ is a set $\{f_N:E_N
\to E'_N\}_{N\geq 1}$ of morphisms of vector bundles such that
$\widetilde{\varphi}_{N} \circ f_{N+1} = f_{N} \circ \widetilde{\varphi}_{N}$.
A morphism of two collections $f:E\to E'$ is called an {\it isomorphism} if it
possesses the inverse morphism.
\medskip

\t\label{Tyurin}{\bf[Barth -- Van de Ven -- Tyurin]} {\it Any
collection $E=\{E_N\}_{N\geq 1}$ of holomorphic vector bundles $E_N$
of a finite rank over $\mathbb{CP}^{N}$ is isomorphic to a direct
sum of collections $E^i=\{E^i_N\}_{N\geq 1}$ of vector bundles
$E^i_N$ of rank $1$.}

\medskip

For collections of rank 2 this result was proved by W.~Barth and A.~Van de Ven
in \cite{BV}, and for collections of an arbitrary finite rank by A.~Tyuirin in
\cite{Tyurin}.

\medskip

The similar question may be considered in the case of complex
supermanifolds. Recall that the {\it projective superspace}
$(M,\mathcal{O})=\mathbb{CP}^{n|m}$ of dimension $n|m$ is a complex
supermanifold with the reduction $M=\mathbb{CP}^{n}$ and the
structure sheaf $\mathcal{O} = \bigwedge\mathcal{L}(-1)^m$, where
$\mathcal{L}(-1)$ is the sheaf of $\mathcal{F}$-modules inverse to
the sheaf $\mathcal{L}(1)$, which corresponds to a hyperplane in
$\mathbb{CP}^n$. The classical homogeneous coordinates $z_0,...,z_n$
on $\mathbb{CP}^{n}$ can be supplemented by odd homogeneous
coordinates $\zeta_1,...,\zeta_m$, giving rise to the system of
homogeneous coordinates on $\mathbb{CP}^{n|m}$.

Let us consider the sequence of projective superspaces
$$
\mathbb{CP}^{1|k_1}\stackrel{\varphi_{1}}{\longrightarrow}
\mathbb{CP}^{2|k_2}\stackrel{\varphi_{2}}{\longrightarrow} \ldots,
$$
where $k_i\le k_{i+1}$ and $\varphi_{i}$ are standard embeddings, i.e any map
$\varphi_{i}: \mathbb{CP}^{i|k_i} \to \mathbb{CP}^{i+1|k_{i+1}}$ is given in
homogeneous coordinates $(z_{j},\, \zeta_r)$ and $(z'_{s},\, \zeta'_t)$ on
$\mathbb{CP}^{i|k_i}$ and $\mathbb{CP}^{i+1|k_{i+1}}$ respectively by
$$
\begin{array}{c}
 z'_{s}= z_{s}, \,\, s=1,\cdots, i, \,\,z_{i+1}=0;\\
\zeta'_t= \zeta_t,\,\, t=1,\cdots, k_i,\,\, \zeta'_t=0,\,\,
t=k_i+1,\cdots, k_{i+1}.
\end{array}
$$

We study collections $\mathcal{E}=\{\mathcal{E}_n\}_{n\geq 1}$ of locally free
sheaves $\mathcal{E}_n$ of a finite rank over $\mathbb{CP}^{n|k_n}$, $n\geq 1$,
such that $\widetilde{\varphi}_{n} (\mathcal{E}_{n+1}) = \mathcal{E}_n$. A
morphism of two collections and their direct sum are defined similarly to the
classical case. We are going to prove the following theorem:

\medskip

\t\label{Tyurin_V} {\it Any collection $\mathcal{E}=\{\mathcal{E}_n\}_{n\geq
1}$ of locally free sheaves $\mathcal{E}_n$ of a finite rank over
$\mathbb{CP}^{n|k_n}$ is isomorphic to a direct sum of collections
$\mathcal{E}^i=\{\mathcal{E}^i_n\}_{n\geq 1}$ of locally free sheaves
$\mathcal{E}^i_n$ of rank $1|0$ or $0|1$.

}

\medskip

\noindent{Proof}. Note that
$\mathcal{E}_{\red}=\{(\mathcal{E}_n)_{\red}\}$ is the collection of
locally free sheaves such that $ \widetilde{(\varphi_i)}_{\red}
((\mathcal{E}_{i+1})_{\red}) = (\mathcal{E}_i)_{\red}$ and
$(\varphi_i)_{\red}: \mathbb{CP}^i \to \mathbb{CP}^{i+1}$ are
standard embeddings. By Theorem \ref{Tyurin} we have
$\mathcal{E}_{\red} \simeq \bigoplus_r \mathcal{S}^r$, where
$\mathcal{S}^r = \{\mathcal{S}^r_n\}$ is a collection of locally
free sheaves of rang $1$ (and of super-rank $1|0$ or $0|1$). Hence
the collection $\gr\mathcal{E}= \{\gr\mathcal{E}_n\}$, where we
identify $\gr\mathcal{E}_n = \mathcal{O}_{\mathbb{CP}^{n}}\otimes
(\mathcal{E}_n)_{\red}$, is isomorphic to the collection $\{
\mathcal{O}_{\mathbb{CP}^{n}}\otimes \bigoplus_r \mathcal{S}^r_n\}$.

Our aim is to show that $\mathcal{E} \simeq \gr\mathcal{E}$. Using
Lemma \ref{Aut_E_(p)=} and the well-known fact:
$H^1(\mathbb{CP}^n,\mathcal L(r))=\{0\}$ for $n>1$ and any $r\in
\mathbb{Z}$, we conclude that
$H^1(\mathbb{CP}^n,\mathcal{E}nd^{\mathcal{O}}_{p}(\gr\mathcal{E}_n))=\{0\}$
for $p\geq 1$ and $n>1$.
 Hence, by
the sequence (\ref{sequense Aut split}) we get
$$
H^1(\mathbb{CP}^n,\mathcal{A}ut^{\mathcal{O}}_{(1)}
(\gr\mathcal{E}_n))=\{\epsilon\}\,\,\, \text{for $n>1$.}
$$
It follows by Theorem \ref{Aut_E_(gr)} that the following
isomorphisms
$$
f_n: \mathcal{E}_n\stackrel{\sim}{\longrightarrow}
 \gr\mathcal{E}_n =
\sum_r \mathcal{O}_{\mathbb{CP}^{n}}\otimes\mathcal{S}^r_n.
$$
exist. Let us show that we can choose the isomorphisms $f_n$ such
that they commute with pullbacks of the bundles. Fix an isomorphism
$f_n$. Let us construct an isomorphism
$$
f'_{n+1}: \mathcal{E}_{n+1}\stackrel{\sim}{\longrightarrow}
\mathcal{O}_{\mathbb{CP}^{n+1}}\otimes (\mathcal{E}_{n+1})_{\red}
$$
 such that $\widetilde{\varphi}_n\circ f'_{n+1} = f_n\circ \widetilde{\varphi}_n$.
 Denote by
 $\mathcal{I}_n$ the sheaf of ideals corresponding to the
 subsupermanifold $\varphi_n: \mathbb{CP}^{n|k_n} \to
 \mathbb{CP}^{n+1|k_{n+1}}$. By definition we have
 $$
 \begin{array}{c}
 \mathcal{E}_{n} = \widetilde{\varphi}_n(\mathcal{E}_{n+1}) = \varphi^*_{\red} (\mathcal{E}_{n+1} / \mathcal{I}_n
 \mathcal{E}_{n+1}), \\
  \gr\mathcal{E}_{n} = \widetilde{\varphi}_n(\gr\mathcal{E}_{n+1}) = \varphi^*_{\red} (\gr\mathcal{E}_{n+1} / \mathcal{I}_n
 \gr\mathcal{E}_{n+1}).
\end{array}
 $$
Denote by $\mathcal B_n$ the sheaf of automorphisms of the sheaf of
$\mathcal{O}_{\mathbb{CP}^{n+1}}/\mathcal{I}_n
\mathcal{O}_{\mathbb{CP}^{n+1}}$-modules
$\gr\mathcal{E}_{n+1}/\mathcal{I}_n \gr\mathcal{E}_{n+1}$ and by
$(\mathcal B_n)_{(1)}$ the subsheaf of $\mathcal B_n$:
$$
(\mathcal B_n)_{(1)} := \{ a\in \mathcal B_n \mid a(v) = v\mod
(\gr\mathcal{E}_{n+1}/\mathcal{I}_n \gr\mathcal{E}_{n+1})_{(1)} \},
$$
 where $(\gr\mathcal{E}_{n+1}/\mathcal{I}_n \gr\mathcal{E}_{n+1})_{(1)}$ is the image of $(\gr\mathcal{E}_{n+1})_{(1)}$ by the natural homomorphism.
 Note that we have $\operatorname{sup}((\mathcal{B}_n)_{(1)}) = \varphi_{\red}(\mathbb{CP}^n)$ and
 $\varphi_{\red}^* ((\mathcal{B}_n)_{(1)}) = \mathcal{A}ut^{\mathcal{O}_{\mathbb{CP}^n}}_{(1)} (\gr\mathcal{E}_{n})$.

 Further, any
 automorphism from $\mathcal{A}ut^{\mathcal{O}_{\mathbb{CP}^n}}_{(1)}
(\gr\mathcal{E}_{n+1})$ preserves $\mathcal{I}_n\gr
 \mathcal{E}_{n+1}$. Hence, we have the map
$$
F_n: \mathcal{A}ut^{\mathcal{O}_{\mathbb{CP}^n}}_{(1)}
(\gr\mathcal{E}_{n+1}) \to (\mathcal{B}_n)_{(1)},
$$
which is surjective as a sheaf homomorphism because we always can find locally preimage of elements from $(\mathcal{B}_n)_{(1)}$.
Denote by $\mathcal{A}_n$ the kernel
of $F_n$.
Let us choose a Stein cover $\mathcal{U}=\{U_i\}$ of
$\mathbb{CP}^{n+1}$ such that
$$
0\to \mathcal{A}_n(U_i) \to
\mathcal{A}ut^{\mathcal{O}_{\mathbb{CP}^{n+1}}}_{(1)}
(\gr\mathcal{E}_{n+1})(U_i) \to (\mathcal{B}_n)_{(1)}(U_i) \to 0.
$$
 is exact for any $i$. Assume also that $\mathcal{U}$ satisfies
 conditions of the proof of Theorem \ref{Aut_E}.
  Denote
by
$$
(g_{ij}^n)\in H^1(\mathcal U, (\mathcal{B}_n)_{(1)})\,\,\,
\text{and}\,\,\, (g_{ij}^{n+1})\in H^1(\mathcal U,
\mathcal{A}ut^{\mathcal{O}_{\mathbb{CP}^n}}_{(1)}
(\gr\mathcal{E}_{n+1}))
$$
the cocycles corresponding to $\mathcal{E}_n$ and
$\mathcal{E}_{n+1}$ by Theorem \ref{Aut_E_(gr)}. Recall that $
g_{ij}^n = \delta_i^n \circ (\delta_j^n)^{-1}$, where $\delta_i^n :
\mathcal{E}_n|_{U_i} \to \gr\mathcal{E}_n|_{U_i}$ is the isomorphism
from Theorem \ref{Aut_E} assuming is addition $\sigma_i = \id$ for
any $i$. Similarly, $g_{ij}^{n+1} = \delta_i^{n+1} \circ
(\delta_j^{n+1})^{-1}$. Since
$\widetilde{\varphi}(\mathcal{E}_{n+1}) = \mathcal{E}_{n}$, we may
assume that $\widetilde{\varphi}^n \circ \delta_i^{n+1}|_{U_i} =
\delta_i^{n} \circ \widetilde{\varphi}^n|_{U_i}$. Therefore,
$F_n(g_{ij}^{n+1}) = g_{ij}^{n}$.

 We have shown that
$(g_{ij}^{n}) \sim \epsilon$ hence there are $\alpha^n_i\in
\mathcal{B}_{(1)}(U_i)$ such that $(\alpha^n_i)^{-1}\circ g_{ij}^{n}
\circ \alpha^n_j = \id$. Using the surjectivity of $F_n|_{U_i}$, we
may choose $\alpha^{n+1}_i \in F_n^{-1}(\alpha^n_i)$. Then
$(h_{ij})\in H^1(\mathcal{U}, \mathcal{A}_n)$, where $h_{ij}=
(\alpha^{n+1}_i)^{-1}\circ g_{ij}^{n+1} \circ \alpha^{n+1}_j$. It is
easy to see that
$$
\begin{array}{rl}
\mathcal{A}_n = \exp(& (\mathcal{I}_n)_{\bar 0} \otimes
((\mathcal{E}_{\red})_{\bar 0} \otimes (\mathcal{E}_{\red})_{\bar
1}^*  \oplus (\mathcal{E}_{\red})_{\bar 1} \otimes
(\mathcal{E}_{\red})_{\bar 0}^*) \oplus\\
&(\mathcal{I}_n)_{\bar 1} \otimes ((\mathcal{E}_{\red})_{\bar 0}
\otimes (\mathcal{E}_{\red})_{\bar 0}^*  \oplus
(\mathcal{E}_{\red})_{\bar 1} \otimes (\mathcal{E}_{\red})_{\bar
1}^* ).
\end{array}
$$
Therefore, we get as for
$\mathcal{A}ut^{\mathcal{O}_{\mathbb{CP}^n}}_{(1)}
(\gr\mathcal{E}_{n})$ that $H^1(\mathbb{CP}^{n+1},
\mathcal{A}_n)=\{\epsilon\}$. Therefore, there are $\beta_i\in
\mathcal{A}_n(U_i)$ such that $h_{ij} = \beta_i\circ \beta_j^{-1}$.
Denote
$$
f'_{n+1}|_{U_i} := \beta_i^{-1}\circ (\alpha_i^{n+1})^{-1} \circ
\delta_i^{n+1}.
$$
By construction, we have $\widetilde{\varphi}_n \circ f'_{n+1} =
f_n\circ \widetilde{\varphi}_n$. The proof is complete.$\Box$

\medskip

\noindent {\it 3.3 About the Birkhoff -- Grothendieck Theorem for
supermanifolds.}

In this subsection we will show that the Birkhoff -- Grothendieck
Theorem:

\smallskip

 \noindent{\it Any finite rank vector bundle on the complex projective space $\mathbb{CP}^1$ is isomorphic to a direct sum of line
bundles},
 \smallskip

 \noindent does not hold true for the projective superspace $\mathbb{CP}^{1|n}$, where $n\geq
1$. Denote by $\mathcal{O}_{n}$ the structure sheaf of
$\mathbb{CP}^{1|n}$ and by $i_n$ the standard embedding
$\mathbb{CP}^{1|1} \to \mathbb{CP}^{1|n}$, $n\geq 1$. Clearly, there
is a map $j_n:\mathbb{CP}^{1|n}\to \mathbb{CP}^{1|1}$, $n\geq 1$,
such that $j_n^*: \mathcal{O}_{1} \to \mathcal{O}_{n}$ is injective
and $j_n\circ i_n=\id$. Let $\mathcal{E}_1$ be a locally free sheaf
of $\mathcal{O}_{1}$-modules. Denote
$$
\mathcal{E}_n:=
\mathcal{O}_{n}\otimes_{j_n^*(\mathcal{O}_{1})}\mathcal{E}_1.
$$
Then $\mathcal{E}_n$ is also locally free and $\mathcal{E}_n$ is an
extension of $\mathcal{E}_1$. In other words, we have proved that
any locally free sheaf on $\mathbb{CP}^{1|1}$ admits an extension to
$\mathbb{CP}^{1|n}$. It follows that to prove our assertion it is
enough to show that there exists a locally free sheaf of
$\mathcal{O}_{1}$-modules of rank $\geq 2$, which is not a direct
sum of two lines bundles.

Let us study firstly line bundles on $\mathbb{CP}^{1|1}$. By
(\ref{sequense Aut split}) we get that
$\mathcal{A}ut^{\mathcal{O}_1}_{(1)} \gr\mathcal{E} \simeq
\mathcal{E}nd^{\mathcal{O}}_{1} \gr\mathcal{E}$ for any rank and
from Lemma \ref{Aut_E_(p)=} it follows that
$\mathcal{E}nd^{\mathcal{O}}_{1} \gr\mathcal{E}=\{0\}$ if $\rank
\gr\mathcal{E}=1|0$ or $0|1$. Therefore, by Theorem \ref{Aut_E_(gr)}
any line bundle $\mathcal{E}$ is isomorphic to $\gr \mathcal{E}$.

Further, let $(\mathcal{E}_{\red})_{\bar 0} = \mathcal{L}(0)$,
$(\mathcal{E}_{\red})_{\bar 1} = \mathcal{L}(-1)$ and
$\mathcal{E}_{\gr} =  \mathcal{O}_{1} \otimes
((\mathcal{E}_{\red})_{\bar 0} \oplus (\mathcal{E}_{\red})_{\bar
1})$. Then
$$
H^1(\mathbb{CP}^{1},\mathcal{E}nd^{\mathcal{O}}_{1}
\mathcal{E}_{\gr})\simeq H^1(\mathbb{CP}^{1}, \mathcal{L}(-2))\simeq
\mathbb{C}.
$$
Using the fact that the unit 1-cohomology class is a fixed point for
the action of $H^0(\mathbb{CP}^{1},
\mathcal{A}ut^{\mathcal{O}_1}_{0} \mathcal{E}_{\gr})$ on
$H^1(\mathbb{CP}^{1}, \mathcal{A}ut^{\mathcal{O}_1}_{(1)}
\mathcal{E}_{\gr})$, we see that there is a locally free sheaf of
$\mathcal{O}_{1}$-modules $\mathcal{E}$ such that $\gr
\mathcal{E}=\mathcal{E}_{\gr}$ but $\mathcal{E}$ is not isomorphic
to $\mathcal{E}_{\gr}$.

\bigskip

\begin{center}
    {\bf 4. The tangent sheaf of a split supermanifold.}
\end{center}

Let us recall some well-known facts about the tangent sheaf
$\mathcal{T}$ of a split supermanifold $(M,\mathcal{O})\simeq (M,
\bigwedge \mathcal{G})$. First, the sheaf $\mathcal{T}$ is
$\mathbb{Z}$-graded (not only $\mathbb{Z}_2$-graded):
$$
\mathcal{T}=\bigoplus_{p\geq -1}\mathcal{T}_p,
$$
where
$$
\mathcal{T}_p:=\{ v\in \mathcal{T} \mid v(\mathcal{O}_q)\subset
\mathcal{O}_{p+q} \,\,\text{for all}\,\, q\geq 0 \},\,\,p\geq -1.
$$
Second, the following sequence
\begin{equation}\label{exact sequence for tangent sheaf}
0\to \bigwedge^{p+1}\mathcal{G} \otimes
\mathcal{G}^*\stackrel{\delta}{\longrightarrow} \mathcal{T}_p
\stackrel{\gamma}{\longrightarrow} \bigwedge^{p}\mathcal{G} \otimes
\Theta\to 0,\,\,p\geq -1,
\end{equation}
where $\Theta$ is the tangent sheaf of $M$, is exact (see
\cite{Oni_Tran Lie}). The mapping $\gamma$ is the restriction of a
derivation of degree $p$ onto the subsheaf $\mathcal{F}\subset
\mathcal{O}$ and $\delta$ identifies any sheaf homomorphism
$\mathcal{G}\to \bigwedge^{p+1}\mathcal{G}$ with a derivation of
degree $p$ that is zero on $\mathcal{F}$.

Denote by $\mathbb{G}$ the vector bundle corresponding to
$\mathcal{G}$. As usual by a {\it (holomorphic) connection} in a
vector bundle $\mathbb{G}\to M$ over a complex manifold $M$, we mean
a bilinear map
$$
\nabla : \Theta \times \mathcal{G} \to \mathcal{G}
$$
satisfying the following conditions:
\begin{itemize}
  \item $\nabla_{fX}s= f\nabla_{X}s$,
  \item $\nabla_{X}(fs) = f \nabla_{X}s + X(f)s$,
\end{itemize}
where $f\in \mathcal{F}$, $X\in \Theta$ and $s\in \mathcal{G}$. If
$\nabla$ and $\nabla'$ are connections in $\mathbb{G}\to M$ and
$\mathbb{G}'\to M$ respectively, the {\it tensor product connection}
$\nabla\otimes \nabla'$ in $\mathbb{G}\otimes \mathbb{G}'$ is well
defined. Recall that
$$
(\nabla\otimes \nabla'_X)(s\otimes s')= \nabla_X(s)\otimes s' +
s\otimes \nabla'_X(s').
$$
It is easy to see that the tensor product connection
$\nabla\otimes\cdots \otimes \nabla$ in $\mathbb{G}\otimes \cdots
\otimes \mathbb{G}$ ($p$-times) induces the {\it wedge product
connection} $\wedge^p \nabla $ in $\bigwedge^p\mathbb{G}$, $p>0$.

Let $\nabla$ be a connection on $\mathbb{G}$. Then to each $X\in
\Theta$ we may assign a vector field $Y_X$ on $(M,\mathcal{O})\simeq
(M,\bigwedge \mathcal{G})$ of degree $0$ defined by
$$
Y_X (f) = X(f), \,\,f\in \mathcal{F},\quad Y_X (f) = \wedge^p
\nabla(f), \,\,f\in \bigwedge^p\mathcal{G},
$$
The Leibniz rule for $Y_X$ follows from the definitions of a
connection and a wedge product connection. Consider the sequence
(\ref{exact sequence for tangent sheaf}) for $p=0$
\begin{equation}\label{exact sequence for tangent sheaf p=0}
0\to \mathcal{G} \otimes
\mathcal{G}^*\stackrel{\delta}{\longrightarrow} \mathcal{T}_0
\stackrel{\gamma}{\longrightarrow}  \Theta\to 0.
\end{equation}
We have just shown that the connection $\nabla$ defines the
splitting of (\ref{exact sequence for tangent sheaf p=0}) by
$X\mapsto Y_X$. The converse statement is also true: if we have a
splitting $i$ of (\ref{exact sequence for tangent sheaf p=0}), we
may define the connection $\nabla_i$ by
$$
(\nabla_i)_X(s):= i(X)(s), \,\,s\in \mathcal{G}.
$$

Note that the curvature tensor of $\nabla= \nabla_i$
$$
R(X,Y)=\nabla_X\circ \nabla_Y - \nabla_Y\circ \nabla_X -
\nabla_{[X,Y]} = ([i(X),i(Y)] - i([X,Y]))|_{\mathcal{G}}
$$
measures the defection of $i$ to be a homomorphism of sheaves of Lie
algebras.

\medskip

\t\label{Prop_Tangent sheaf} {\it Let $(M,\mathcal{O}_M) \simeq (M,
\bigwedge \mathcal{G})$ be a (holomorphic) split supermanifold and
$\mathcal{T}$ the tangent sheaf. The following conditions are
equivalent:
\begin{enumerate}
  \item the sheaf $\mathcal{T}$ corresponds to the unit 1-cohomology class with values in
  $\mathcal{A}ut^{\mathcal{O}}_{(1)}
\gr\mathcal{T}$ by the Theorem \ref{Aut_E_(gr)};
  \item the sequence (\ref{exact sequence for tangent sheaf p=0})
splits;
  \item $\mathcal{G}$ possesses a (holomorphic) connection.
\end{enumerate}

}

\medskip

\noindent{Proof}. By the discussion above we have to prove only that
$\mathcal{T}$ corresponds to the trivial 1-cocycles of
$H^1(M,\mathcal{A}ut^{\mathcal{O}}_{(1)} \gr\mathcal{T})$ if and
only if the sequence (\ref{exact sequence for tangent sheaf p=0})
splits. Let $\theta_0: \Theta \to \mathcal{T}_0$ be a splitting of
(\ref{exact sequence for tangent sheaf p=0}). Then the sequence
(\ref{exact sequence for tangent sheaf}) splits for all $p\geq 0$,
we may define the splitting $\theta_p: \bigwedge^{p}\mathcal{G}
\otimes\Theta \to \mathcal{T}_p$ by $\theta_p(f\otimes v)= f
\theta_0(v)$. It follows that
$$
\mathcal{T}_p\simeq \bigwedge^{p}\mathcal{G} \otimes\Theta +
\bigwedge^{p+1}\mathcal{G} \otimes \mathcal{G}^*.
$$
 Hence,
$$
\mathcal{T} \simeq \bigwedge\mathcal{G} \otimes
(\mathcal{G}^*+\Theta) \simeq \bigwedge\mathcal{G} \otimes
(\mathcal{T}_{\red})=\gr \mathcal{T}.
$$

Conversely, since the unit cocycle of
$H^1(M,\mathcal{A}ut^{\mathcal{O}}_{(1)} \gr\mathcal{T})$ is a fixed
point with respect to the action of $H^0(M,
\mathcal{A}ut^{\mathcal{O}}_{0} \gr\mathcal{T})$, there is an
isomorphism $\Phi: \mathcal{T}\to \gr \mathcal{T}$ such that $\gr
\Phi = \id$ (see proof of Theorem \ref{Aut_E}). It follows that the
following diagram is commutative
$$
\begin{CD}
\mathcal{T}_{\bar 0}@>{\Phi|_{\mathcal{T}_{\bar 0}}}>> (\gr \mathcal{T})_{\bar 0}\\
@V{\pi}VV @VV{\pr}V\\
\mathcal{T}_{\bar 0}/(\mathcal{J}\mathcal{T})_{\bar 0}@=
\mathcal{T}_{\bar 0}/(\mathcal{J}\mathcal{T})_{\bar 0}
\end{CD},
$$
where $\pr$ is the projection of
$$
\gr\mathcal{T}= \bigoplus_{p\geq 0}
(\mathcal{J}^{p}\mathcal{T})_{\bar
0}/(\mathcal{J}^{p+1}\mathcal{T})_{\bar 0}+ \bigoplus_{p\geq 0}
(\mathcal{J}^{p}\mathcal{T})_{\bar
1}/(\mathcal{J}^{p+1}\mathcal{T})_{\bar 1}
$$
onto $\mathcal{T}_{\bar 0}/(\mathcal{J}\mathcal{T})_{\bar 0}$ and
$\pi$ is the natural projection. Further, by the definitions of all
morphisms the following diagram is also commutative
$$
\begin{CD}
\mathcal{T}_{\bar 0}@>{\pi}>> \mathcal{T}_{\bar 0}/(\mathcal{J}\mathcal{T})_{\bar 0}\\
@V{\pr_{\mathcal{T}_{0}}}VV @VV{\tau}V\\
\mathcal{T}_{0}@>{\gamma}>> \Theta
\end{CD},
$$
where $\tau$ is an isomorphism defined by $v+
(\mathcal{J}\mathcal{T})_{\bar 0} \mapsto \pr_{\mathcal{F}}\circ
v|_{\mathcal{F}}$. Denote by $i$ the natural embedding
$\mathcal{T}_{\bar 0}/(\mathcal{J}\mathcal{T})_{\bar 0}
\hookrightarrow (\gr \mathcal{T})_{\bar 0}$. We may define a
splitting of (\ref{exact sequence for tangent sheaf p=0}) by
$\pr_{\mathcal{T}_{0}}\circ (\Phi|_{\mathcal{T}_{\bar 0}})^{-1}
\circ i\circ \tau^{-1}$. The proof is complete.$\Box$

\bigskip

\begin{center}
    {\bf 5.  A spectral sequence}
\end{center}

An important problem is to calculate the cohomology group
$H^*(M,\mathcal{E})$ of a locally free sheaf of
$\mathcal{O}$-modules $\mathcal{E}$ on a supermanifold
$(M,\mathcal{O})$. If $(M,\mathcal{O})$ is split, then $\mathcal{E}$
is a locally free sheaf of $\mathcal{F}$-modules on $M$, and its
cohomology group can be calculated in many cases using the well
elaborated tools of complex analytic geometry. In non-split case
these methods cannot be applied directly, but we can use the
associated split supermanifold $(M,\gr \mathcal{O})$ and the sheaf
$\gr \mathcal{E}$.

\medskip

\noindent {\it 5.1 Quasi-derivations.}

Let $(M,\mathcal{O})$ be an arbitrary supermanifold and
$\mathcal{E}$ a locally free sheaf on $(M,\mathcal{O})$. Let us take
an even vector field $\Gamma\in \mathcal{T}_{\bar 0}(U)$ on a
superdomain $(U,\mathcal{O}|U)\subset(M,\mathcal{O})$. A
$\mathbb{Z}_2$-graded vector spaces sheaf homomorphism
$A_{\Gamma}:\mathcal{E}|U \to \mathcal{E}|U$ is called a {\it
$\Gamma$-derivation} if $A_{\Gamma}(fv)=\Gamma(f)v+f A_{\Gamma}(v)$,
$f\in \mathcal{O}|U$ and $v\in \mathcal{E}|U$. A homomorphism of
$\mathbb{Z}_2$-graded sheaf of vector spaces $B:\mathcal{E} \to
\mathcal{E}$ will be called a {\it quasi-derivation} if it is a
$\Gamma$-derivation for a certain $\Gamma$. Denote by
$\mathcal{QD}er\mathcal{E}$ the sheaf of quasi-derivations. It is a
sheaf of Lie algebras with respect to the commutator $[A_{\Gamma},
B_{\Upsilon}]:= A_{\Gamma}\circ B_{\Upsilon} - B_{\Upsilon} \circ
A_{\Gamma}$. The sheaf $\mathcal{QD}er\mathcal{E}$ possesses the
double filtration:
$$
\begin{array}{rl}
\mathcal{QD}er_{(p)(q)}\mathcal{E}:= \{ A_{\Gamma}\in
\mathcal{QD}er\mathcal{E} \mid& A_{\Gamma}(\mathcal{E}_{(r)})\subset
\mathcal{E}_{(r+p)},\,\,
\Gamma(\mathcal{J}^s)\subset \mathcal{J}^{s+q} \\
& \text{for all}\,\, \, r,s\in \mathbb{Z}\}.
\end{array}
$$
The map
$$
\exp: \mathcal{QD}er_{(1)(2)}\mathcal{E} \to \mathcal{QA}ut_{(1)(2)} \mathcal{E}
$$
is an isomorphism of sheaves of sets.
Let us consider the subsheaf $\mathcal{QD}er_{k, k}\gr\mathcal{E}$ of
$\mathcal{QD}er_{(k)(k)}\gr\mathcal{E}$ defined by
$$
\begin{array}{rl}
\mathcal{QD}er_{k, k}\gr\mathcal{E} := &\{ A_{\Gamma}\in
\mathcal{QD}er_{(k)(k)}\gr\mathcal{E} \mid
A_{\Gamma}(\gr\mathcal{E}_{r})\subset \gr\mathcal{E}_{r+k},\,\,\\
&\rule{0pt}{5mm}\Gamma(\gr\mathcal{O}_s)\subset
\gr\mathcal{O}_{s+k}\,\,\, \text{for all}\,\, \, r,s\in
\mathbb{Z}\}.
\end{array}
$$
Note that $\mathcal{QD}er_{k, k}\gr\mathcal{E} =
\mathcal{E}nd^{\gr\mathcal{O}}_k\gr \mathcal{E}$ if $k$ is odd.

Denote by $\mu_k$, $k\geq 1$, the following mapping:
$$
\mu_k: \mathcal{QA}ut_{(k)(2)} \gr\mathcal{E} \to \mathcal{QD}er_{k,
k}\gr\mathcal{E},
$$
$$
\mu_k(a_{\gamma}) = \bigoplus_q \pr_{q+k} \circ A_{\Gamma} \circ
\pr_{q},
$$
where $a_{\gamma} = \exp (A_{\Gamma})$ and $\pr_{k}: \gr\mathcal{E}
\to \gr\mathcal{E}_k$ is the natural projection. The kernel of this
map is $\mathcal{QA}ut_{(k+1)(2)} \gr\mathcal{E}$. Moreover, the
following sequence
$$
0\to \mathcal{QA}ut_{(k+1)(2)} \gr\mathcal{E} \longrightarrow
\mathcal{QA}ut_{(k)(2)} \gr\mathcal{E}
\stackrel{\mu_k}{\longrightarrow} \mathcal{QD}er_{k, k}\gr\mathcal{E}\to 0
$$
is exact.
Denoting by $H_{(k)}(\gr\mathcal{E})$ the image of the natural
mapping $$
H^1(M,\mathcal{QA}ut_{(k)(2)}\gr \mathcal{E})\to
H^1(\mathcal{QA}ut_{(1)(2)}\gr \mathcal{E}),
$$ we get the filtration:
$$
H^1(M,\mathcal{QA}ut_{(1)(2)}\gr \mathcal{E}) = H_{(1)}(\gr\mathcal{E})
\supset H_{(2)}(\gr\mathcal{E})\supset \ldots .
$$
Take $a_{\gamma}\in H_{(1)}(\gr\mathcal{E})$. We define the order of
$a_{\gamma}$ the maximal one of the numbers $k$ such that
$a_{\gamma}\in H_{(k)}(\gr\mathcal{E})$. The {\it order} of a
locally free sheaf $\mathcal{E}$ of $\mathcal{O}$-modules on a
supermanifold $(M,\mathcal{O}_M)$ is by definition the order of the
corresponding cohomology class.

\medskip

\noindent {\it 5.2 The spectral sequence.}

Let $\mathcal{E}$ be a vector superbundle on a supermanifold
$(M,\mathcal{O})$ of dimension $n|m$. Now we will construct a
spectral sequence for the cohomology of the sheaf $\mathcal{E}$. We
fix an open Stein cover $\frak U = (U_i)_{i\in I}$ of $M$ and
consider the corresponding \v{C}ech cochain complex $C^*(\frak
U,\mathcal{E}) = \bigoplus_{p\ge 0} C^p(\frak U,\mathcal{E})$.

The $\Bbb Z_2$-grading of $\mathcal{E}$ gives rise to the $\Bbb
Z_2$-gradings in $C^*(\frak U,\mathcal{E})$ and $H^*(M,\mathcal{E})$
given by
\begin{equation}\label{grading of H}
\aligned C_{\bar 0}(\frak U,\mathcal{E}) &= \bigoplus_{q\ge 0}
C^{2q}(\frak U,\mathcal{E}_{\bar 0})\oplus\bigoplus_{q\ge 0}
C^{2q+1}(\frak U,\mathcal{E}_{\bar 1}),\\
C_{\bar 1}(\frak U,\mathcal{E}) &= \bigoplus_{q\ge 0} C^{2q}(\frak
U,\mathcal{E}_{\bar 1})\oplus\bigoplus_{q\ge 0}
C^{2q+1}(\frak U,\mathcal{E}_{\bar 0}).\\
H_{\bar 0}(M,\mathcal{E}) &= \bigoplus_{q\ge 0}
H^{2q}(M,\mathcal{E}_{\bar 0})\oplus
\bigoplus_{q\ge 0} H^{2q+1}(M,\mathcal{E}_{\bar 1}),\\
H_{\bar 1}(M,\mathcal{E}) &= \bigoplus_{q\ge 0}
H^{2q}(M,\mathcal{E}_{\bar 1})\oplus \bigoplus_{q\ge 0}
H^{2q+1}(M,\mathcal{E}_{\bar 0}).
\endaligned
\end{equation}

The filtration (\ref{filtr_rassloenie}) for $\mathcal{E}$ gives rise
to the filtration
\begin{equation}\label{filtr ckomplex}
C^*(\frak U,\mathcal{E}) = C_{(0)}\supset\ldots\supset C_{(p)}
\supset\ldots\supset C_{(m+1)} = 0
\end{equation}
of this complex by the subcomplexes
$$
C_{(p)} = C^*(\frak U,\mathcal{E}_{(p)}).
$$
Denoting by $H(M,\mathcal{E})_{(p)}$ the image of the natural
mapping $H^*(M,\mathcal{E}_{(p)})\to H^*(M,\mathcal{E})$, we get the
filtration
\begin{equation}\label{filtr H}
H^*(M,\mathcal{E}) = H(M,\mathcal{E})_{(0)}\supset\ldots \supset
H(M,\mathcal{E})_{(p)}\supset \ldots.
\end{equation}
Denote by $\gr H^*(M,\mathcal{E})$ the bigraded group associated
with the filtration (\ref{filtr H}); its bigrading is given by
$$
\gr H^*(M,\mathcal{E}) = \bigoplus_{ p,q\ge 0}\gr_p
H^q(M,\mathcal{E}).
$$
By the general procedure, invented by Leray, the filtration
(\ref{filtr ckomplex}) gives rise to a spectral sequence of bigraded
groups $E_r$ converging to $E_{\infty}\simeq \gr
H^*(M,\mathcal{E})$. It is constructed in the following way.

For any $p,r\ge 0$, define the vector spaces
$$
C^p_r = \{c\in C_{(p)}\,|\,dc\in C_{(p+r)}\}.
$$
Then, for a fixed $p$, we have
$$
C_{(p)} = C^p_0\supset\ldots\supset C^p_r\supset
C^p_{r+1}\supset\ldots.
$$
The $r$-th term of the spectral sequence is defined by
$$
E_r = \bigoplus_{p=0}^m E^p_r,\;r\ge 0,
$$
where
$$
E^p_r = C^p_r/C^{p+1}_{r-1} + dC^{p-r+1}_{r-1}.
$$
Since $d(C^p_r)\subset C^{p+r}_r$, $d$ induces a derivation $d_r$ of
$E_r$ of degree $r$ such that $d_r^2 = 0$. Then $E_{r+1}$ is
naturally isomorphic to the homology algebra $H(E_r,d_r)$. Denoting
$Z_r = \Ker d_r$, we have the natural mapping $\kappa^r_{r+1}:
Z_r\to E^{r+1}$. For any $s > r$, denote $\kappa^r_s =
\kappa^{s-1}_s\circ\ldots\circ\kappa^r_{r+1}$ (this composition is
not defined on the entire $Z_r$).

The $\Bbb Z_2$-grading (\ref{grading of H}) in $C^*(\frak
U,\mathcal{E})$ gives rise to certain $\Bbb Z_2$-gradings in $C^p_r$
and $E^p_r$, turning $E_r$ into a superspace. Clearly, the
coboundary operator $d$ in $C^*(\frak U,\mathcal{E})$ is odd. It
follows that the coboundary $d_r$ is odd for any $r\ge 0$.

The superspaces $E_r$ are also endowed with a second $\Bbb
Z$-grading. Namely, for any $q\in \Bbb Z$, set
$$
\aligned
C^{p,q}_r &= C^p_r\cap C^{p+q}(\frak U,\mathcal{E}),\\
E^{p,q}_r &= C^{p,q}_r/C^{p+1,q-1}_{r-1} + dC^{p-r+1,q+r-2}_{r-1}.
\endaligned
$$
Then
$$
E_r = \bigoplus_{p,q} E^{p,q}_r.
$$
Clearly,
\begin{equation}\label{d_r(E^(p,q)_r)subset}
d_r(E^{p,q}_r)\subset E^{p+r,q-r+1}_r
\end{equation}
for any $r,\,p,\,q$.

One sees easily that $C^{p,q}_r = 0$ for all $p$ and $r$ if $q\le
-(m+1)$. Therefore, for a fixed $q$, we have $d(C^{p,q}_r) = 0$ for
all $r\ge q+m+2$. This implies that $\kappa^r_{r+1}: E^{p,q}_r\to
E^{p,q}_{r+1}$ is an isomorphism for all $p$ and $r\ge r_0(q) =
q+m+2$. Setting $E^{p,q}_{\infty} = E^{p,q}_{r_0(q)}$, we get the
bigraded superspace
$$
E_{\infty} = \bigoplus_{p,q} E^{p,q}_{\infty}.
$$

Now we prove certain properties of the spectral sequence $(E_r)$.
Some of them are well known and are valid in a more general
situation.

\medskip

\l\label{E_0,E_1 and E_2} {\it The first two terms of the spectral
sequence $(E_r)$ can be identified with the following bigraded
spaces:
$$
\aligned
E_0 &= C^*(\frak U,\gr\mathcal{E}),\\
E_1 &= H^*(M,\gr\mathcal{E}).
\endaligned
$$
Here
$$
\aligned
E_0^{p,q} &= C^{p+q}(\frak U,(\gr\mathcal{E})_p),\\
E_1^{p,q} &=  H^{p+q}(M,(\gr\mathcal{E})_p).
\endaligned
$$

}

\noindent{Proof.} By definition, we have
$$
E_0^p = C_{(p)}/C_{(p+1)},\; p\ge 0,
$$
where the coboundary operator $d_0$ of degree 0 is induced by $d:
C_{(p)}\to C_{(p)}$. On the other hand, the exact sequence
$$
0\to\mathcal{E}_{(p+1)}\to\mathcal{E}_{(p)}\to \gr\mathcal{E}_p\to 0
$$
and Theorem B for Stein supermanifolds imply the exact sequence
$$
0\to\mathcal{E}_{(p+1)}(U)\to\mathcal{E}_{(p)}(U)\to
\gr\mathcal{E}_p(U)\to 0
$$
for any Stein open subset $U\subset M$. Therefore
$$
C^*(\frak U,(\gr\mathcal{E})_p)\simeq C_{(p)}/C_{(p+1)} = E_0^p,\;
p\ge 0.
$$
One sees easily that this is an isomorphism of complexes and that
the resulting isomorphism $C^*(\frak U,\gr\mathcal{E})\simeq E_0$ is
an isomorphism of bigraded spaces. It follows that
$$
E_1\simeq H(E_0,d_0)\simeq H^*(\frak U,\gr\mathcal{E})\simeq
H^*(M,\gr\mathcal{E}).\Box
$$

\medskip

\medskip

\l\label{E_infinity} {\it There is the following identification of
bigraded algebras:
$$
E_{\infty} = \gr H^*(M,\mathcal{E}),
$$
where
$$
E_{\infty}^{p,q} = \gr_p H^{p+q}(M,\mathcal{E}).
$$

}

\noindent{Proof.} Clearly, for $r\ge r_0(q)$ we have $C_r^{p,q} =
Z^{p+q}(\frak U, \mathcal{E}_{(p)})$. It follows that
$$
\aligned E^{p,q}_{\infty} &= Z^{p+q}(\frak U,\mathcal{E}_{(p)})/
Z^{p+q}(\frak U,\mathcal{E}_{(p+1)}) + dC^{p+q-1}(\frak
U,\mathcal{E})\cap
Z^{p+q}(\frak U,\mathcal{E}_{(p)})\\
&= H^{p+q}(M,\mathcal{E})_{(p)}/(Z^{p+q}(\frak
U,\mathcal{E}_{(p+1)})/
dC^{p+q-1}(\frak U,\mathcal{E})\cap Z^{p+q}(\frak U,\mathcal{E}_{(p+1)})\\
&= H^{p+q}(M,\mathcal{E})_{(p)}/H^{p+q}(M,\mathcal{E})_{(p+1)} =
\gr_p H^{p+q}(M,\mathcal{E}).\Box
\endaligned
$$

\medskip

\noindent{\bf Corollary.} {\it If $M$ is compact, then
$$
\dim H^k(M,\mathcal{E}) = \sum_{p+q = k}\dim  E_{\infty}^{p,q}.
$$

}

\noindent{Proof.} In fact, if $M$ is compact, then all cohomology
groups with values in a coherent analytic sheaf on $(M,\mathcal{O})$
or $M$ are of finite dimension.$\Box$

\medskip

Now we prove our main result concerning the first non-zero
coboundary operators among $d_1,\;d_2,\ldots$. We may suppose that
for each $i\in I$ there exists an isomorphism of sheaves $\sigma_i:
\mathcal{O}|U_i\to\gr \mathcal{O}|U_i$, inducing the identity
isomorphism $\gr \mathcal{O}|U_i\to\gr\mathcal{O}|U_i$.

By Theorem \ref{Aut_E}, a locally free sheaf of
$\mathcal{O}$-modules $\mathcal{E}\to (M,\mathcal{O})$ corresponds
to the cohomology class $a_{\gamma}$ of the $1$-cocycle
$((a_{\gamma})_{ij})\in Z^1(\frak U, \mathcal{QA}ut_{(1)(2)}
\gr\mathcal{E})$, where $(a_{\gamma})_{ij} =
\delta_i\circ\delta_j^{-1}$. If the order of $(a_{\gamma})_{ij}$ is
equal to $k$, then we may choose $\delta_i,\; i\in I$, in such a way
that $((a_{\gamma})_{ij})\in Z^1(\frak U,\mathcal{QA}ut_{(k)(2)}
\gr\mathcal{E})$. We can write $a_{\gamma} = \exp A_{\Gamma}$, where
$A_{\Gamma}\in C^1(\frak U,\mathcal{QD}er_{(1)(2)} \gr\mathcal{E})$.

We will identify the differential spaces $(E_0,d_0)$ and $(C^*(\frak
U, \gr\mathcal{E}),d)$ via the isomorphism of Proposition
\ref{E_0,E_1 and E_2}. Clearly, $\delta_i: \mathcal{E}_{(p)}|U_i\to
\gr\mathcal{E}_{(p)}|U_i = \sum_{r\ge p}\gr\mathcal{E}_r|U_i$ is an
isomorphism of sheaves for any $i\in I,\; p\ge 0$. These local sheaf
isomorphisms permit us to define an isomorphism of graded cochain
groups
$$
\psi: C^*(\frak U,\mathcal{E})\to C^*(\frak U,\gr\mathcal{E})
$$
such that
$$
\psi: C^*(\frak U,\mathcal{E}_{(p)})\to C^*(\frak
U,\gr\mathcal{E}_{(p)}),\;p\ge 0.
$$
We give it by
$$
\psi(c)_{i_0\ldots i_q} = \delta_{i_0}( c_{i_0\ldots i_q})
$$
for any $(i_0,\ldots,i_q)$ such that $U_{i_0}\cap\ldots\cap
U_{i_q}\ne \emptyset$. In general, $\psi$ is not an isomorphism of
complexes. Nevertheless, we can express explicitly the coboundary
$d$ of the complex $C^*(\frak U,\mathcal{E})$ by means of $d_0$ and
$a_{\gamma}$.

\medskip

\l\label{psi d psi^(-1)c} {\it For any $c\in C^q(\frak
U,\gr\mathcal{E}) = \oplus_p E_0^{q-p,p}$, we have
$$
(\psi(d\psi^{-1}(c)))_{i_0\ldots i_{q+1}} = (d_0c)_{i_0\ldots
i_{q+1}} + ((a_{\gamma})_{i_0i_1} -\id )(c_{i_1\ldots i_{q+1}}).
$$

}

\noindent{Proof.} We can write
$$
\aligned (d\psi^{-1}(c))_{i_0\ldots i_{q+1}} &= \sum_{\alpha
=0}^{q+1} (-1)^{\alpha}
\psi^{-1}(c)_{i_0\ldots \hat i_{\alpha}\ldots i_{q+1}}\\
&= \sum_{\alpha =1}^{q+1}(-1)^{\alpha}\psi^{-1}(c)_{i_0\ldots \hat
i_{\alpha}\ldots i_{q+1}}
+ \psi^{-1}(c)_{i_1\ldots i_{q+1}}\\
&= \delta_{i_0}^{-1}(\sum_{\alpha =1}^{q+1} (-1)^{\alpha}
c_{i_0\ldots \hat i_{\alpha}\ldots i_{q+1}}) +
\delta_{i_1}^{-1}(c_{i_1\ldots i_{q+1}}) \\
&= \delta_{i_0}^{-1}((d_0c)_{i_0\ldots i_{q+1}} - c_{i_1\ldots
i_{q+1}}) + \delta_{i_1}^{-1}( c_{i_1\ldots i_{q+1}}).
\endaligned
$$
Therefore
$$
\aligned (\psi(d\psi^{-1}(c)))_{i_0\ldots i_{q+1}} &= \delta_{i_0}
(d\psi^{-1}(c))_{i_0\ldots i_{q+1}}\\
&= (d_0c)_{i_0\ldots i_{q+1}} - c_{i_1\ldots i_{q+1}}
+ (a_\gamma)_{i_0i_1} (c_{i_1\ldots i_{q+1}}) \\
&= (d_0c)_{i_0\ldots i_{q+1}} + ((a_\gamma)_{i_0i_1} - \id)
(c_{i_1\ldots i_{q+1}}).
\endaligned
$$
This implies our assertion.$\Box$

\medskip

This proposition makes it possible to calculate the spectral
sequence $(E_r)$ whenever $d_0$ and the cochain $a_{\gamma}$ are
known. Now we find the explicit form of certain coboundary operators
$d_r,\; r\ge 1$.

\medskip

\t\label{Teor_d} {\it Suppose that the locally free sheaf of
$\mathcal{O}$-modules $\mathcal{E}\to (M,\mathcal{O}_M)$ has order
$k$ and denote by $a_{\gamma}$ the cohomology class corresponding to
$\mathcal{E}$ by Theorem \ref{Aut_E}. Then $d_r = 0$ for $r =
1,\ldots,k-1$, and $d_{k} = \mu_k(a_{\gamma})$.

}

\medskip

\noindent{Proof.} Take a cocycle $c\in E^{p,q-p}_0,\; d_0c = 0$, and
denote by $c^*$ its cohomology class in $E^{p,q-p}_1$. Clearly, $c$
and $c^*$ are represented by the cochain $\psi^{-1}(c)\in C^p_0$. By
Proposition \ref{psi d psi^(-1)c},
$$
(\psi(d\psi^{-1}(c)))_{i_0\ldots i_{q+1}} = ((a_\gamma)_{i_0i_1} -
\id) (c_{i_1\ldots i_{q+1}}).
$$
Now we see that
$$
(\psi(d\psi^{-1}(c)))_{i_0\ldots i_{q+1}} =
\mu_k(a_{\gamma})_{i_0i_1}( c_{i_1\ldots i_{q+1}}) + u_{i_0\ldots
i_{q+1}},
$$
where $u\in C_{(p+k+1)}$. This means that
$$
\psi(d\psi^{-1}(c)) = \mu_k(a_{\gamma})( c) + u,
$$
whence $d_1 = d_2 =\ldots = d_{(k-1)} = 0$. Identifying $E_{k}$ with
$E_1$, we also see that $d_{k}c^*$ is represented by the cochain
$\psi^{-1}(\mu_k(a_{\gamma})( c))$. It follows that
$$
d_{2k}c^* = \mu_k(a_{\gamma})(c^*).\Box
$$

\smallskip

\noindent{\it Arkady Onishchik}

\noindent {\textsc{Yaroslavl University\\
150 000 Yaroslavl, Russia}}\\
\noindent {\emph{E-mail address:} \verb"aonishch@aha.ru",

\bigskip

\noindent{\it Elizaveta Vishnyakova}

\noindent {\textsc{Max-Planck-Institut f\"{u}r Mathematik\\
P.O.Box: 7280\\
53072 Bonn\\
Germany}}

 \noindent {\emph{E-mail address:}
\verb"VishnyakovaE@googlemail.com", \verb"Liza@mpim-bonn.mpg.de"}

\end{document}